\documentclass[11pt]{article}

\usepackage[margin=1in]{geometry}
\usepackage{amsmath}
\usepackage{amsthm}
\usepackage{color}
\usepackage{graphicx}
\usepackage{newtxtext}
\usepackage[subscriptcorrection]{newtxmath}
\usepackage{nicefrac}
\usepackage{comment}
\usepackage{tikz}
\usepackage[round,authoryear]{natbib}
\usepackage[hidelinks]{hyperref}
\graphicspath{{./art/}}

\newtheorem{theorem}{Theorem}
\newtheorem{proposition}{Proposition}
\theoremstyle{definition}
\newtheorem{definition}{Definition}

\def\T{{ \mathrm{\scriptscriptstyle T} }}

\DeclareMathOperator{\var}{var}
\DeclareMathOperator{\cov}{cov}

\title{Whittle--Mat\'ern metrics and isotropic Gaussian processes on compact metric graphs}
\author{David Bolin\quad and\quad Alexandre B. Simas\\[6pt]
\small Statistics Program, King Abdullah University of Science and Technology (KAUST),\\
\small Thuwal 23955, Saudi Arabia\\
\small \texttt{david.bolin@kaust.edu.sa, alexandre.simas@kaust.edu.sa}}
\date{}

\begin{document}
\maketitle

\begin{abstract}
Isotropic Gaussian processes on networks have a covariance that is a function of a distance between locations, but valid constructions have so far required graphs with Euclidean edges and the resistance metric. We introduce a one-parameter family of metrics $d_\kappa$, $\kappa\ge0$, on arbitrary compact metric graphs, including those with loops and parallel edges, defined as the variograms of first-order Whittle--Mat\'ern fields, where $\kappa=0$ corresponds to an intrinsic limiting field. We show that these Whittle--Mat\'ern metrics are invariant under insertion and removal of vertices of degree two, that they decrease monotonically in $\kappa$, and that they can be evaluated by sparse linear solves. The metric $d_0$ coincides with the resistance metric on graphs with Euclidean edges, and thus extends it to general compact metric graphs without any modification of the graph. For every completely monotone function $r$ and every $\kappa\ge0$, we show that $r(d_\kappa)$ is a positive semidefinite function. This yields isotropic Gaussian processes on general compact metric graphs. For $\kappa>0$ the metric levels off beyond the screening length $\kappa^{-1}$, which can be used to induce long-range dependence in  short-range covariance models, and $\kappa$ can be estimated from data together with the parameters of $r$. The models are compared on daily traffic maps from the street network of Paris and the motorway network around London. For both, $\kappa$ is estimated to be positive, $d_\kappa$ increases the maximised log-likelihood relative to $d_0$ for every covariance function, and it gives the best predictions at randomly held-out locations.
\end{abstract}

\section{Introduction}\label{sec:intro}

Many spatial data sets are observed on networks rather than in Euclidean space. Water quality and temperature are measured along river networks \citep{Hoef2006,peterson2013}, and traffic speeds, accidents and other events are recorded on street networks \citep{ROAD2,ROAD3,baddeley2021review,moller2022lgcp}. For such data, distances should be measured along the network, and the observations and the locations where predictions are needed lie at arbitrary positions on the roads or rivers, not only at the junctions. The natural domain is therefore a compact metric graph $\Gamma$, i.e., a finite set of vertices $\mathcal{V}$ joined by a finite set of edges $\mathcal{E}$, where each edge $e$ is a curve of length $\ell_e$, every point of every edge is a location, and distance is measured by the geodesic metric $d(\cdot,\cdot)$, the length of the shortest path along the edges. On such a domain the vertex set is largely a choice of representation. Declaring a sensor location to be a vertex inserts a vertex of degree two and changes the graph representation, even though the resulting network is canonically isometric to the original one. Real networks typically also contain parallel edges, such as two roads joining the same pair of junctions, and loops. A statistical model should therefore assign the same joint distribution to a given set of locations irrespective of which other points have been declared vertices, and it should not require the network to be simplified before it can be used.

The models of choice in Euclidean space are isotropic Gaussian processes, whose covariance $r\{\hat{d}(s,t)\}$ depends on the locations only through their distance. They have constant variance, their parameters are easily interpreted, and many parametric families are available. Transferring them to networks is, however, not straightforward, because a function $r$ which gives a valid covariance in Euclidean space need not do so when it is composed with a distance on a graph. For the geodesic metric, \citet{anderes2020isotropic} showed that the exponential covariance fails to be positive semidefinite, for some value of its range parameter, on any graph with Euclidean edges in which two points are joined by three paths with pairwise disjoint interiors. 
They therefore introduced the resistance metric, which extends the effective resistance of electrical network theory \citep{klein_randic1993} from the vertices to the points on the edges, and showed that every completely monotone function of this metric is a valid covariance function on graphs with Euclidean edges. The other main approaches give up isotropy. The convolution models of \citet{Hoef2006} are designed for stream networks and are defined on metric trees. The Whittle--Mat\'ern fields of \citet{BSW2022}, defined through a stochastic partial differential equation on the graph, exist on every compact metric graph, can be differentiable, and have Markov properties \citep{BSW_Markov} which make them computationally attractive, but their covariance is not a function of a distance and their marginal variance varies over the network.

Isotropic models are thus available only under a restriction on the network. The requirement of Euclidean edges means that there is at most one edge between each pair of vertices, that there are no loops, and that the length of every edge equals the geodesic distance between its endpoints. It is satisfied for linear networks and for metric trees, but not for general metric graphs such as complex street networks \citep{ROAD2,ROAD3}. A key reason for the requirement is that the resistance metric of \citet{anderes2020isotropic} is constructed from the Laplacian matrix of the underlying simple graph, and it is therefore not defined on a metric graph that violates these conditions. Such a graph can be brought into the required form by inserting vertices of degree two on the offending edges, but the metric is then only available through a modified representation of the network, and the inserted vertices, which serve no other purpose, enlarge the matrices that enter the computations. The construction also provides a single geometry, with no parameter that controls how the network structure is translated into distance.

In this work we take a different route, which starts from the Whittle--Mat\'ern fields and applies to every compact metric graph. For $\kappa>0$, let $u_\kappa$ denote a first-order Whittle--Mat\'ern field, with covariance operator $(\kappa^2-\Delta_\Gamma)^{-1}$, where $\Delta_\Gamma$ is the Kirchhoff--Laplacian \citep{berkolaiko_kuchment2013}. The values of $u_\kappa$ at the vertices have a sparse precision matrix $Q^\kappa$, which is available in closed form and is consistent under insertion and removal of vertices of degree two 
\citep{BSW2022,BSW2023_AOS}. We define the \emph{Whittle--Mat\'ern metrics} as the variograms
\begin{equation}\label{eq:metric}
d_\kappa(s,t)=\var\{u_\kappa(s)-u_\kappa(t)\},\qquad s,t\in\Gamma,\quad \kappa\ge0,
\end{equation}
where $d_0$ is obtained as the limit as $\kappa\downarrow0$ of $d_{\kappa}$.
We show that each $d_\kappa$ is a metric on $\Gamma$ which is invariant under insertion and removal of vertices of degree two, and that $d_\kappa(s,t)$ is continuous and strictly decreasing in $\kappa$ for $s\neq t$. The metrics require no condition on the graph, and they are evaluated directly from the given representation by sparse linear solves with $Q^\kappa$, where the only vertices that need to be added are the locations of interest.

The case $\kappa=0$ connects the construction to the graph Laplacian. As $\kappa\downarrow0$, the matrix $Q^\kappa$ converges to the matrix $Q$ whose elements are determined by the conductances $1/\ell_e$. We refer to $Q$ as the generalised graph Laplacian of $\Gamma$, since it is defined for every representation of the metric graph, including those with loops and parallel edges. For a simple graph it is the Laplacian matrix of the electrical network in which every edge is a resistor with resistance $\ell_e$ \citep{doyle_snell1984}, which is also used in the theory of metrized graphs \citep{baker_faber2006,baker_rumely2007,cinkir2011}, and the elimination of vertices from such matrices is known as Kron reduction \citep{dorfler_bullo2013}. The matrices $Q$ associated with different degree-two subdivisions are exactly consistent under Schur-complement elimination of the inserted vertices. Thus the construction is compatible with metrically equivalent graph representations rather than being tied to a particular choice of vertex set. We link this construction to the Whittle--Mat\'ern fields as the limit of their precision matrices at the vertices. As a consequence, $d_0$ is the effective resistance between arbitrary points of $\Gamma$, which satisfies $d_0\le d$ with equality for all pairs of points if and only if $\Gamma$ is a tree. On graphs with Euclidean edges, $d_0$ coincides with the resistance metric by \citet{anderes2020isotropic}. Thus $d_0$ extends the resistance metric to arbitrary compact metric graphs without first replacing the given representation by one with Euclidean edges, and identifies it as one member of a family of metrics generated by the Whittle--Mat\'ern fields. For $\kappa>0$, the reciprocal $\kappa^{-1}$ acts as a screening length: the metric no longer grows with the separation between the locations but levels off beyond $\kappa^{-1}$. On the real line, for example, $d_\kappa(s,t)=\kappa^{-1}(1-e^{-\kappa|s-t|})$, whereas $d_0(s,t)=|s-t|$.

The Whittle--Mat\'ern metrics yield isotropic Gaussian processes on general compact metric graphs. We show that
\begin{equation}\label{eq:kernel_family}
\rho(s,t)=r\{d_\kappa(s,t)\},\qquad s,t\in\Gamma,
\end{equation}
is positive semidefinite on every compact metric graph, for every $\kappa\ge0$ and every completely monotone function $r$, which includes the exponential, the Mat\'ern with smoothness $\nu\le\nicefrac12$, and the generalised Cauchy families. This generalises the result of \citet{anderes2020isotropic} in two directions: it removes the restriction to graphs with Euclidean edges, and it replaces the single resistance geometry by the family $d_\kappa$. The resulting Gaussian processes have constant variance $r(0)$, and the parameters of $r$ determine how dependence decays with dissimilarity, as in Euclidean space, whereas $\kappa$ determines the dependence between distant locations. Since $d_\kappa$ levels off beyond the screening length, the correlation between locations that are far apart does not decay to zero but remains at a positive level controlled by $\kappa$. Thus $\kappa>0$ induces long-range dependence, in the sense of a positive lower bound on correlations, in models which are otherwise short-range,  and the short-range model in the resistance metric is recovered as $\kappa\downarrow0$. In Euclidean space, long-range dependence typically requires an additional latent component or a different class of covariance functions, whereas here it is obtained through a single additional parameter, which can be estimated from data together with the parameters of $r$ by maximum likelihood. Since the metrics are invariant under subdivision, the model does not depend on which locations have been declared vertices, so estimation at the observation locations and prediction at further locations refer to the same process. This is in contrast to covariances defined as matrix functions of a vertex Laplacian, such as the diffusion and graph Mat\'ern kernels used in learning on graphs \citep{kondor_lafferty2002,smola_kondor2003,borovitskiy2021}, which are defined only at the vertices and change on the existing vertices when a vertex is inserted. 



\section{A generalised graph Laplacian}\label{sec:graph_laplacian}

A connected compact metric graph $\Gamma$ is obtained from a finite connected multigraph $G=(\mathcal V,\mathcal E)$ with vertex set $\mathcal{V}=\{v_1,\ldots,v_n\}$ and edge set $\mathcal{E}$, where parallel edges and self-loops are allowed, by associating with each edge $e\in\mathcal E$ a length $\ell_e>0$ and a copy of the interval $[0,\ell_e]$, and identifying interval endpoints according to the incidence relations of $G$. The graph is equipped with the path metric induced by the edge lengths. Thus every interior point of an edge, and not only the vertices, is a location of $\Gamma$. Inserting a degree-two vertex in the interior of an edge changes the graph representation, but the resulting geometric realization is canonically isometric to the original one as the subdivision preserves edge lengths, points of the network, and all geodesic distances. We refer to two representations related in this way as degree-two subdivisions of one another. For locations $s_1,\ldots,s_m\in\Gamma$, we write $\Gamma(s_1,\ldots,s_m)$ for the degree-two subdivision of $\Gamma$ obtained by inserting those locations that are not already vertices, and $\mathcal V^\star$ for its vertex set.

The Laplacian matrix of $\Gamma$ associated with this representation, which we refer to as the \emph{generalised graph Laplacian} of $\Gamma$, is the $n\times n$ matrix $Q$ with elements
\begin{equation}\label{eq:Q0}
Q_{ij} = \begin{cases}
\sum\limits_{e \in \mathcal{E}_{v_i}} \dfrac{1}{\ell_e}\mathbb{I}(\bar{e}\neq \underline{e}) & i=j,\\[4pt]
\sum\limits_{e \in \mathcal{E}_{v_i} \cap \mathcal{E}_{v_j}} -\dfrac{1}{\ell_e} & i\neq j,
\end{cases}
\end{equation}
where $\mathcal{E}_{v}$ denotes the set of graph edges incident to $v$, each edge appearing once, and $\mathbb{I}(\bar e\neq\underline e)$ equals $1$ precisely for non-loop edges. Distinct parallel edges are distinct members of $\mathcal E_v$. For differential vertex conditions we write $\mathcal H_v$ for the multiset of edge ends incident to $v$: a non-loop edge contributes one end at each endpoint, while a self-loop based at $v$ contributes two ends.

The form of \eqref{eq:Q0} is motivated by the fact that $Q$ arises as the $\kappa\downarrow 0$ limit of the precision matrix of a Gaussian Whittle--Mat\'ern field on $\Gamma$, which gives $Q$ a natural probabilistic interpretation. Recall from \citet{BSW2022} that the Gaussian Whittle--Mat\'ern field with parameters $\kappa>0$ and $\alpha>\nicefrac12$ is the centred Gaussian process on $\Gamma$ with covariance operator $\mathcal{C}=(\kappa^2-\Delta_\Gamma)^{-\alpha}$, where $\Delta_\Gamma$ is the Kirchhoff--Laplacian, which acts as the second derivative on each edge with domain
  \begin{equation}\label{eq:Kirchoff}
\mathcal{D}(\Delta_{\Gamma}) = \left\{ f \in \widetilde{H}^2(\Gamma)\cap C(\Gamma) : \sum_{h \in \mathcal{H}_v} \partial_h f(v) = 0 \text{ for every }v\in\mathcal V \right\}.
\end{equation} 
Here $\widetilde{H}^2(\Gamma)$ is the space of functions on $\Gamma$ that are twice weakly differentiable on the edges, $C(\Gamma)$ is the space of  continuous functions on $\Gamma$, and $\partial_h$ denotes the outward derivative at the edge end $h$. Counting edge ends is essential for self-loops, whose two endpoint traces both enter the Kirchhoff condition. We focus on the special case $\alpha=1$, in which the Whittle--Mat\'ern field $u_\kappa$ is a first-order Markov process \citep{BSW_Markov}: by \citet{BSW2023_AOS}, the values of $u_\kappa$ at the vertices of $\Gamma$ form a centred Gaussian vector with $n\times n$ precision matrix $Q^\kappa$ given by
\begin{equation}\label{eq:expprec}
Q^\kappa_{ij} = 2\kappa\begin{cases}
\displaystyle\sum_{e \in \mathcal{E}_{v_i}} \left[\left(\dfrac{1}{2} + \dfrac{e^{-2\kappa \ell_e}}{1-e^{-2\kappa \ell_e}}\right) \mathbb{I}(\bar{e}\neq \underline{e})+ \tanh(\kappa\ell_e/2)\, \mathbb{I}(\bar{e}=\underline{e})\right] & i=j,\\[6pt]
\displaystyle\sum_{e \in \mathcal{E}_{v_i} \cap \mathcal{E}_{v_j}} -\dfrac{e^{-\kappa \ell_e}}{1-e^{-2\kappa \ell_e}} & i\neq j,
\end{cases}
\end{equation}
The field $u_\kappa$ is unchanged by the addition or removal of vertices of degree $2$, and the precision matrix of the original vertices is recovered from that of a subdivision by Schur complement elimination of the inserted vertices \citep{BSW2022}. 
%
%
As stated in Proposition~\ref{prop:Q_properties}\,(vi) below, the precision matrix \eqref{eq:expprec} converges entry-wise to $Q$ in \eqref{eq:Q0} as $\kappa\downarrow 0$, and $Q$ is the Dirichlet-to-Neumann matrix associated with the energy of $-\Delta_\Gamma$; see Proposition~\ref{prop:edge_DN}. This ties the matrix $Q$ directly to the Kirchhoff--Laplacian on the metric graph.

The term generalised graph Laplacian refers to its metric-graph interpretation, in which it is defined for every representation of $\Gamma$, rather than to a new algebraic notion of multigraph Laplacian. For a simple weighted graph $G=(\mathcal{V},\mathcal{E})$ with positive edge weights $\{w_{ij}\}$, the weighted combinatorial graph Laplacian $L^G$ is
\begin{equation}\label{eq:standard_laplacian}
L^G_{ij} = \begin{cases}
\sum_{k:(v_i,v_k)\in\mathcal{E}} w_{ik}, & i=j,\\
-w_{ij}, & (v_i,v_j)\in\mathcal{E},\\
0, & \text{otherwise}.
\end{cases}
\end{equation} 
With the conductance weights $w_{ij}=1/\ell_e$, this is the Laplacian matrix of the electrical network in which every edge is a resistor with resistance $\ell_e$ \citep{doyle_snell1984}, and it is the matrix associated with a metric graph in the theory of metrized graphs, where the vertex set is chosen so that the graph is simple \citep{baker_faber2006,cinkir2011}. Weighted multigraph Laplacians are also standard: with unit vertex weights, \citet[Definition~3.6 and Lemma~3.8(ii)--(iii)]{AliniaeifardWangVanWilligenburg2024} show that loops do not contribute and that parallel edges combine by summing their weights, and \eqref{eq:Q0} is the matrix of the weighted graph that is associated with a metric graph with loops and parallel edges in the literature on quantum graphs \citep{kostenko_nicolussi2023}. What is specific here is the relation of $Q$ to the continuous metric graph: its weights are fixed by the edge lengths, it is the limit of the precision matrices of the Whittle--Mat\'ern fields, and the resulting family of matrices is compatible with subdivision. The results cited above concern either the matrix alone, simple representations, or single edges, and we are not aware of a statement of these properties for general representations with loops and parallel edges, which is the form needed in the following sections. They are given in the next proposition, with proofs in the appendix.

\begin{proposition}\label{prop:Q_properties}
Let $\Gamma$ and $Q$ be as above. Then the following hold.
\begin{enumerate}
\item[i.] $Q$ is symmetric and its quadratic form is
\begin{equation}\label{eq:dirichlet_form}
x^\T Q x = \sum_{e=(v_i,v_j)\in\mathcal{E}: v_i\neq v_j}\frac{(x_i-x_j)^2}{\ell_e},\qquad x\in\mathbb{R}^{n}.
\end{equation}
\item[ii.] $Q$ is positive semi-definite, and $\ker(Q) = \mathrm{span}\{\mathbf{1}\}$, where $\mathbf{1}=(1,\ldots,1)^\T\in\mathbb{R}^n$. In particular, $Q$ admits a spectral decomposition with non-negative eigenvalues $0=\mu_0<\mu_1\leq\cdots\leq\mu_{n-1}$.
\item[iii.] If the underlying combinatorial graph of $\Gamma$ is simple, then $Q$ coincides with the standard weighted graph Laplacian $L^G$ in \eqref{eq:standard_laplacian} with conductance weights $w_{ij}=1/\ell_e$ for any edge $e$ with endpoints $v_i$ and $v_j$.
\item[iv.] If $k$ parallel edges $e_1,\ldots,e_k\in\mathcal{E}_{v_i}\cap\mathcal{E}_{v_j}$ connect distinct vertices $v_i$ and $v_j$, then
\begin{equation}\label{eq:parallel_conductance}
Q_{ij} = -\sum_{m=1}^k \frac{1}{\ell_{e_m}},
\end{equation}
which is the parallel-conductance combination of the individual edge conductances $1/\ell_{e_m}$.
\item[v.] Loops contribute zero to all entries of $Q$.
\item[vi.] The matrix $Q$ in \eqref{eq:Q0} is the entry-wise limit of the precision matrix $Q^\kappa$ in \eqref{eq:expprec} as $\kappa\downarrow 0$.
 \item[vii.] Let $\widetilde\Gamma$ be obtained from $\Gamma$ by inserting finitely many degree-two vertices, let $\widetilde Q$ be the corresponding matrix, and order the refined vertex set as the original vertices $\mathcal V$ followed by the inserted vertices $\mathcal I$. Then $\widetilde Q_{\mathcal I\mathcal I}$ is positive definite and
\begin{equation}\label{eq:kron_subdivision}
Q=\widetilde Q_{\mathcal V\mathcal V}-\widetilde Q_{\mathcal V\mathcal I}\widetilde Q_{\mathcal I\mathcal I}^{-1}\widetilde Q_{\mathcal I\mathcal V}.
\end{equation}
\end{enumerate}
\end{proposition}

Parts (i)--(v) are the counterparts, for general representations, of standard properties of weighted Laplacian matrices; see \citet{baker_faber2006} and \citet{cinkir2011} for simple representations of metric graphs. Part (vi) gives $Q$ a probabilistic interpretation as the limit of the Whittle--Mat\'ern precision matrices, including the contributions of parallel edges and loops. Part (vii) is the consistency property needed below: if the vertex representation is refined by inserting interior locations, the original matrix is recovered exactly by Schur complementation. It is related to Kron reduction in electrical network theory, where elimination of vertices by Schur complementation is known to map Laplacian matrices to Laplacian matrices and to preserve the effective resistances between the remaining vertices \citep{dorfler_bullo2013}. The additional content of (vii) is that the reduced matrix is the generalised graph Laplacian of the coarser, canonically isometric representation. Thus the matrices $Q$ associated with different degree-two subdivisions are compatible under Schur complementation rather than being attached to a particular choice of vertex set. Loops, although invisible on a vertex set containing only their base vertex, become ordinary pairs of parallel edge segments as soon as an interior point of the loop is inserted.

\section{Whittle--Mat\'ern metrics}\label{sec:metrics}

\subsection{Definition and properties}\label{sec:WMmetric}

For $\kappa>0$, let $u_\kappa$ be the first-order Gaussian Whittle--Mat\'ern field introduced by \citet{BSW2022}, obtained as the solution of
\begin{equation}\label{eq:wm_spde_alpha1}
(\kappa^2-\Delta_\Gamma)^{\nicefrac{1}{2}}u_\kappa=\mathcal{W},
\end{equation}
where $\mathcal{W}$ denotes Gaussian white noise on $\Gamma$. We define the positive-$\kappa$ Whittle--Mat\'ern metric by the increment variance of $u_\kappa$, and define the endpoint $\kappa=0$ by continuity from positive $\kappa$.
\begin{definition}\label{def:wm_metric}
For $\kappa>0$, the Whittle--Mat\'ern metric $d_\kappa$ is defined by \eqref{eq:metric}. At $\kappa=0$, define
\begin{equation}\label{eq:d0_limit_definition}
 d_0(s,t)=\lim_{\kappa\downarrow0}d_\kappa(s,t).
\end{equation}
\end{definition}

The following theorem clarifies the properties of this metric. In particular, it shows that the limiting definition of $d_0$ is well posed.

\begin{theorem}\label{thm:metric_properties}
Let $\Gamma$ be a compact metric graph. Then the limit in \eqref{eq:d0_limit_definition} exists for every $s,t\in\Gamma$, and with $d_0$ defined by this limit the following hold.
\begin{enumerate}
\item[i.] $d_\kappa(\cdot,\cdot)$ is a metric on $\Gamma$ for every $\kappa\ge0$ and is invariant under degree-two subdivision, under the canonical identification of the two geometric realizations.
\item[ii.] For fixed $s,t\in\Gamma$, $d_\kappa(s,t)$ is continuous in $\kappa\in[0,\infty)$. If $s\neq t$ and $0\le\kappa_1<\kappa_2$, then $d_{\kappa_1}(s,t)>d_{\kappa_2}(s,t)$, and $d_\kappa(s,t)\to0$ as $\kappa\to\infty$.
\item[iii.] $d_0(s,t)\le d(s,t)$ for every $s,t\in\Gamma$, with equality for each $s,t$ if and only if  $\Gamma$ is a metric tree.
\item[iv.] For every $\kappa>0$ and $s,t\in\Gamma$,
\begin{equation}\label{eq:screening_bound}
 d_\kappa(s,t)\le \frac{2}{\kappa}\tanh\!\left\{\frac{\kappa d(s,t)}{2}\right\}<\frac{2}{\kappa}.
\end{equation}
\end{enumerate}
\end{theorem}

For $u_\kappa$ in \eqref{eq:wm_spde_alpha1}, $\kappa$ is a range parameter where larger values correspond to shorter correlation ranges of $u_\kappa$. For a covariance $r(d_\kappa)$ with $r$ completely monotone and non-constant, the effect is the opposite, since $d_\kappa$ decreases in $\kappa$, so that all correlations of the resulting field increase with $\kappa$. We therefore refer to $\kappa^{-1}$ as the screening length of the metric.
For $\kappa>0$, the effect of increasing network separation is attenuated once distances are large relative to $\kappa^{-1}$, whereas $d_0$ has no corresponding screening parameter. Thus a short-range radial covariance composed with $d_\kappa$ can retain non-negligible dependence between locations that are far apart on the network. In this sense, the positive-$\kappa$ metrics allow a form of long-range dependence controlled by $\kappa$. This is discussed in more detail in Section~\ref{sec:isotropic}. 

The distinction can be seen particularly clearly by considering intervals whose length tends to infinity. For $T>0$, let $\Gamma_T=[-T,T]$ be regarded as a compact metric graph, write $d_{\kappa,T}$ for its Whittle--Mat\'ern metric, and fix $s,t\in\mathbb R$ with $T>\max(|s|,|t|)$. A direct calculation, given in the appendix, shows that
\begin{equation}\label{eq:line_metric_limit}
\lim_{T\to\infty}d_{\kappa,T}(s,t)=
\begin{cases}
|s-t|,&\kappa=0,\\[2mm]
\displaystyle \kappa^{-1}\{1-e^{-\kappa|s-t|}\},&\kappa>0.
\end{cases}
\end{equation}
We write $d_{\kappa,\mathbb R}$ for the limiting metrics in \eqref{eq:line_metric_limit}. Then $d_{0,\mathbb R}(s,t)=|s-t|$ grows without bound as $|s-t|\to\infty$, whereas $d_{\kappa,\mathbb R}(s,t)$ is bounded by $\kappa^{-1}$ and converges to $\kappa^{-1}$ as $|s-t|\to\infty$ for every $\kappa>0$. The two limits also have a familiar probabilistic interpretation: $d_{0,\mathbb R}$ is the variogram of a two-sided standard Brownian motion, while $d_{\kappa,\mathbb R}$ is the variogram of the stationary Ornstein--Uhlenbeck process with covariance $(2\kappa)^{-1}e^{-\kappa|s-t|}$, equivalently the first-order Whittle--Mat\'ern field on $\mathbb R$ under the normalization used here.

The precision matrix $Q^\kappa$ in \eqref{eq:expprec} describes $u_\kappa$ only at the vertices. To evaluate $d_\kappa(s_1,s_2)$ for two locations $s_1,s_2\in\Gamma$, we therefore work on the subdivision $\Gamma(s_1,s_2)$.
\begin{proposition}\label{prop:dkappa_computation}
Let $s_1,s_2\in\Gamma$, let $\mathcal V^\star$ be the vertex set of $\Gamma(s_1,s_2)$, and let $e=e_{s_1s_2}\in\mathbb R^{|\mathcal V^\star|}$ have entry $+1$ at $s_1$, entry $-1$ at $s_2$, and zero elsewhere. Let $Q^\kappa$ be the precision matrix in \eqref{eq:expprec} for $\kappa>0$, and set $Q^0=Q$. Then
\begin{equation}\label{eq:dkappa_formula}
d_\kappa(s_1,s_2)=
\begin{cases}
e^\T(Q^\kappa)^{-1}e,&\kappa>0,\\[2mm]
e^\T Q^\dagger e,&\kappa=0.
\end{cases}
\end{equation}
For $\kappa>0$, the equation $Q^\kappa y=e$ has a unique solution and $d_\kappa=e^\T y$. For $\kappa=0$, $Qy=e$ is solvable, any two solutions differ by an additive constant, and the same identity $d_0=e^\T y$ holds for every solution. Thus the metric can be evaluated by sparse linear solves without explicitly forming an inverse or pseudoinverse.
\end{proposition}

\begin{figure}
    \includegraphics[width=\linewidth]{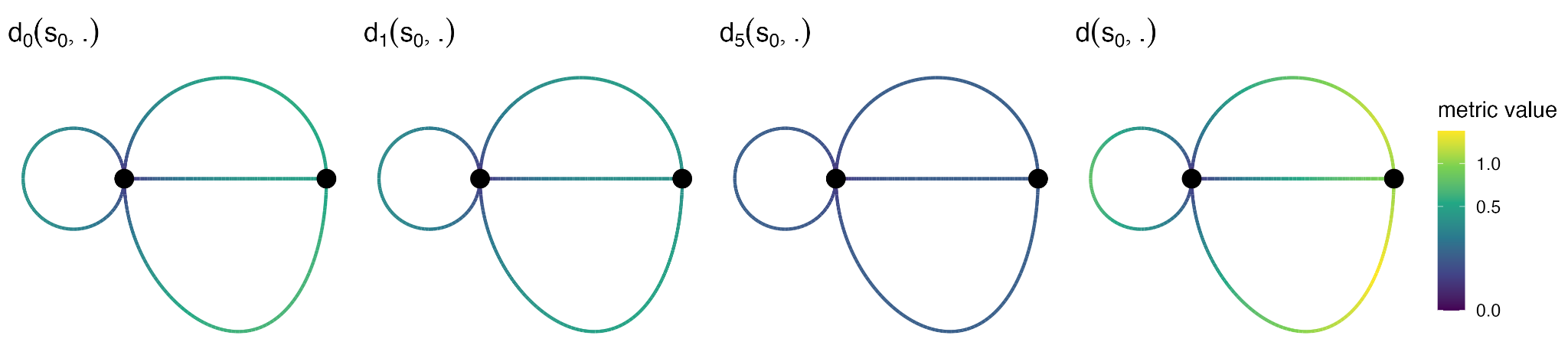}
    \caption{The Whittle--Mat\'ern metrics $d_\kappa(s_0,\cdot)$ for $\kappa\in\{0,1,5\}$ and the geodesic distance $d(s_0,\cdot)$, evaluated on a multigraph with three parallel edges between two vertices and a self-loop. The reference point $s_0$ is the left vertex.}
    \label{fig:metric_multigraph}
\end{figure}
 An example of the metric is shown in Figure~\ref{fig:metric_multigraph}, computed on a metric graph $\Gamma$ with two vertices $v_0=(0,0)$ and $v_1=(1,0)$ connected by three edges, namely a straight segment of length $1$, a semicircular arc of length $\pi/2$, and a longer arc of length approximately $1.85$, together with a self-loop attached at $v_0$. The combinatorial structure of $\Gamma$ contains both parallel edges between $v_0$ and $v_1$ and a loop at $v_0$, so that $\Gamma$ does not have Euclidean edges and falls outside the framework of \citet{anderes2020isotropic}. The four panels of Figure~\ref{fig:metric_multigraph} show the Whittle--Mat\'ern metrics $d_\kappa(s_0,\cdot)$, with $s_0=v_0$, for $\kappa\in\{0,1,5\}$, evaluated by means of Proposition~\ref{prop:dkappa_computation} on a fine mesh of $\Gamma$, together with the geodesic distance $d(s_0,\cdot)$ for comparison. As asserted by Theorem~\ref{thm:metric_properties}, the metrics $d_\kappa$ are well-defined on this multigraph and decrease monotonically in $\kappa$. The geodesic distance, which uses only the shortest path between any two locations, is everywhere larger than $d_0$, with the gap most pronounced on the lower arc and on the loop, where the parallel-conductance effect described in Proposition~\ref{prop:Q_properties}\,(iv) reduces $d_0$ relative to the geodesic.

\subsection{Connection to the resistance metric}\label{sec:resistance_metric}
Recall that a connected compact metric graph has Euclidean edges if it has no loops, at most one edge between each pair of vertices, and $\ell_e=d(u,v)$ for every edge $e$ joining vertices $u$ and $v$. For such a graph, \citet{anderes2020isotropic} defined the resistance metric as the variogram
\begin{equation}\label{eq:res_metric_def}
d_R(s_1,s_2)=\var\{Z_\Gamma(s_1)-Z_\Gamma(s_2)\},\qquad s_1,s_2\in\Gamma,
\end{equation}
of the auxiliary Gaussian process
\begin{equation}\label{eq:aux_process_resmetric}
Z_\Gamma(s)=Z_\mu(s)+\sum_{e\in\mathcal E}Z_e(s),\qquad s\in\Gamma,
\end{equation}
whose components are independent. Here $Z_e$ is a Brownian bridge on the edge $e$, vanishing at its endpoints and outside $e$, and $Z_\mu$ is the linear interpolation along each edge of a centred Gaussian vector at the vertices with precision matrix $L=L^G+e_{v_0}e_{v_0}^\T$, where $L^G$ is the graph Laplacian \eqref{eq:standard_laplacian} with weights $w_{ij}=1/d(v_i,v_j)$, $v_0$ is an arbitrarily chosen vertex, and $e_{v_0}$ is the corresponding unit vector. On a graph with Euclidean edges, $L^G=Q$ by Proposition~\ref{prop:Q_properties}\,(iii), and the following theorem shows that the resistance metric is the Whittle--Mat\'ern metric with $\kappa=0$.

\begin{theorem}\label{thm:resistance_equiv}
For a compact metric graph $\Gamma$ with Euclidean edges, 
$d_0 =d_R$.
\end{theorem}
Thus $d_0$, extends the resistance metric to general compact metric graphs. On a representation without Euclidean edges, the construction of \citet{anderes2020isotropic} is not available, since $L^G$ is defined for simple graphs and uses the weights $1/d(v_i,v_j)$, which differ from the conductances $1/\ell_e$ when $\ell_e>d(v_i,v_j)$, whereas $Q$ is assembled from the actual edges. A representation with Euclidean edges can always be obtained by inserting vertices of degree two, and by Theorem~\ref{thm:metric_properties}\,(i) and Theorem~\ref{thm:resistance_equiv} the resistance metric of any such subdivision coincides with $d_0$, which is thus obtained directly from the original representation.

\begin{figure}[t]
\centering
\begin{tikzpicture}[font=\footnotesize,
  vtx/.style={circle,fill,inner sep=1.3pt},
  ins/.style={circle,draw,fill=white,inner sep=1.3pt}]
\begin{scope}
\node[vtx,label=left:$v_1$] (a1) at (0,0) {};
\node[vtx,label=above left:$v_0$] (a0) at (1.3,0) {};
\node[vtx,label=right:$v_2$] (a2) at (3.5,0) {};
\draw (a1) -- node[below] {$1$} (a0);
\draw (a0) to[bend left=50] node[above] {$\ell$} (a2);
\draw (a0) to[bend right=50] node[below] {$\ell$} (a2);
\end{scope}
\begin{scope}[xshift=6cm]
\node[vtx,label=left:$v_1$] (b1) at (0,0) {};
\node[vtx,label=above left:$v_0$] (b0) at (1.3,0) {};
\node[vtx,label=right:$v_2$] (b2) at (3.5,0) {};
\draw (b1) -- node[below] {$1$} (b0);
\draw (b0) to[bend left=50] node[pos=0.2,above left=-2pt] {$\ell/2$} node[pos=0.8,above right=-2pt] {$\ell/2$} node[ins,pos=0.5] {} (b2);
\draw (b0) to[bend right=50] node[pos=0.2,below left=-2pt] {$\ell/2$} node[pos=0.8,below right=-2pt] {$\ell/2$} node[ins,pos=0.5] {} (b2);
\end{scope}
\end{tikzpicture}
\caption{A metric graph and the subdivision obtained by inserting a vertex at the midpoint of each parallel edge.}
\label{fig:parallel_edges_example}
\end{figure}
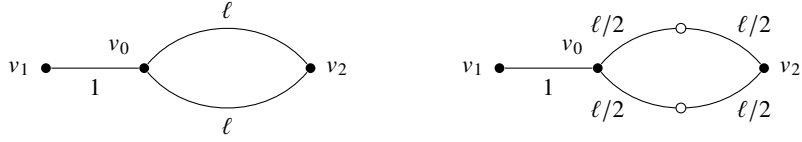

As an example, consider the graph in Figure~\ref{fig:parallel_edges_example}. The conductances of the parallel edges add, so that
\begin{equation}\label{eq:parallel_example_Q}
Q=
\begin{pmatrix}
1+2\ell^{-1}&-1&-2\ell^{-1}\\
-1&1&0\\
-2\ell^{-1}&0&2\ell^{-1}
\end{pmatrix},
\qquad
d_0(v_0,v_2)=\ell/2,\qquad d_0(v_1,v_2)=1+\ell/2,
\end{equation}
whereas the simple graph in which the two parallel edges are replaced by a single edge of length $\ell$ gives the distances $\ell$ and $1+\ell$. The graph can therefore not be replaced by its simple graph on the same vertex set. Inserting a vertex at the midpoint of each parallel edge gives a subdivision with Euclidean edges, and eliminating the two inserted vertices from its Laplacian by \eqref{eq:kron_subdivision} recovers the matrix $Q$ in \eqref{eq:parallel_example_Q}.

\section{Isotropic Gaussian processes}\label{sec:isoGP}

\subsection{Definition and properties}\label{sec:isotropic}
Having defined the Whittle--Mat\'ern metric, we can now use this to define Gaussian processes which are isotropic in this metric. For this, we need the notion of completely   monotone functions.

\begin{definition}
A function $f:[0,\infty)\to\mathbb{R}$ is completely monotone if it is continuous, infinitely differentiable on $(0,\infty)$, and $(-1)^m f^{(m)}(x)\geq0$ for all $x>0$ and $m=0,1,\ldots$
\end{definition}

Our main result is now that any completely monotone function induces a valid isotropic covariance function. The result applies to the full family $d_\kappa$, including the limiting metric $d_0$.

\begin{theorem}\label{thm:cm_positive_def}
    Let $r : [0,\infty)\to \mathbb{R}$ be a completely monotone function, then for any $\kappa\geq 0$, the function $\rho : \Gamma \times \Gamma \to \mathbb{R}$ defined by $\rho(s,t) = r(d_{\kappa}(s,t))$ is a positive semidefinite kernel on $\Gamma\times\Gamma$.
\end{theorem}

\begin{figure}
    \includegraphics[width=\linewidth]{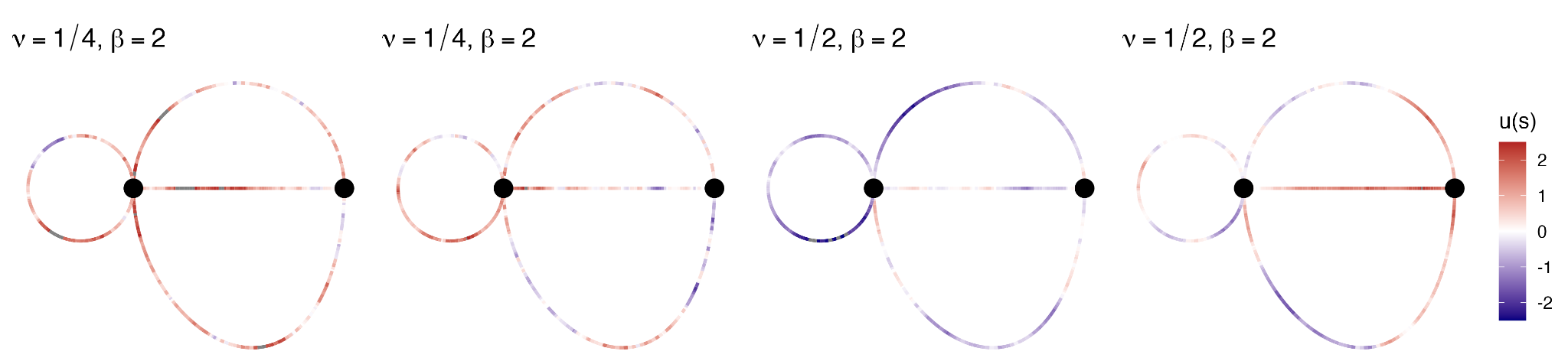}
    \caption{Two independent realisations of a centred Gaussian field with covariance $r(d_0(s,t))$, for $r$ the Mat\'ern covariance function with $\sigma^2=1$, $\beta=2$, and smoothness $\nu=\nicefrac14$ (left two panels) and $\nu=\nicefrac12$ (right two panels). }
    \label{fig:samples_isotropic}
\end{figure}

This generalises the result of \citet{anderes2020isotropic} to the full family $d_\kappa$ and to general compact metric graphs, that may not  have Euclidean edges. Examples of completely monotone functions are given by \citet{moller2022lgcp}; a notable one is the Mat\'ern covariance function
\begin{equation}\label{eq:matern}
r(h) = \frac{\sigma^2}{\Gamma(\nu)\,2^{\nu - 1}}(\beta h)^{\nu} K_{\nu}(\beta h),
\end{equation}
with $\sigma,\beta>0$ and $\nu \in (0,\nicefrac12]$. Figure~\ref{fig:samples_isotropic} shows two independent realisations of a centred Gaussian process with covariance $r\{d_0(s,t)\}$ on the metric graph of Figure~\ref{fig:metric_multigraph}, for this covariance function with  $\nu=\nicefrac14$ and $\nu=\nicefrac12$. As expected, the realisations are rougher for the smaller value of $\nu$.

The covariance $r\{d_\kappa(s,t)\}$ is defined for all $s,t\in\Gamma$ and, by Theorem~\ref{thm:metric_properties}\,(i), does not depend on the representation of $\Gamma$. As for the models of \citet{anderes2020isotropic} on graphs with Euclidean edges, estimation and prediction therefore refer to the same Gaussian process: if the model is fitted at locations $s_1,\ldots,s_n$ and predictions are then required at further locations, these can be inserted as vertices of degree two without changing the covariances among $s_1,\ldots,s_n$. The difference is that this holds for every representation of the network, including those with loops and parallel edges, so that a representation with Euclidean edges does not have to be constructed first.

In a covariance $r\{d_\kappa(\cdot,\cdot)\}$, the parameter $\kappa$ enters alongside the parameters of $r$. To make its role transparent, consider the exponential model
\begin{equation}\label{eq:expmodel_dkappa}
\cov\{X(s),X(t)\}=\sigma^2\exp\{-\beta d_\kappa(s,t)\},
\qquad s,t\in\Gamma,
\end{equation}
where $\sigma^2$ is the marginal variance and $\beta$ is the range parameter of the radial covariance. The real-line limit \eqref{eq:line_metric_limit} satisfies $d_{\kappa,\mathbb R}(s,t)\sim |s-t|$ as $t\to s$ for every $\kappa\ge0$, so $\kappa$ does not alter the short-distance behaviour of the model. Its effect is instead on the dependence at large separations: as $|s-t|\to\infty$, the correlation $\exp\{-\beta d_{\kappa,\mathbb R}(s,t)\}$ tends to zero for $\kappa=0$ but to $\exp(-\beta/\kappa)>0$ for $\kappa>0$. More generally, if $r(h)\to0$ as $h\to\infty$, the limiting correlation is $r(1/\kappa)/r(0)>0$ for $\kappa>0$. On a general compact metric graph, Theorem~\ref{thm:metric_properties}\,(iv) gives $d_\kappa<2/\kappa$, so every correlation in the model is at least $r(2/\kappa)/r(0)$, irrespective of the size of the network. Thus $\kappa>0$ induces long-range dependence through the metric, even though the radial covariance is short-range as a function of its argument. The parameter $\kappa$ therefore provides a parsimonious way of controlling this form of long-range dependence, and it can be estimated jointly with the parameters of $r$. The dependence is of a specific type. If $r(h)=\int e^{-\lambda h}\mu(d\lambda)$, then on the real line $r\{d_{\kappa,\mathbb R}(s,t)\}=\sum_{n\ge0}w_n e^{-n\kappa|s-t|}$ with $w_n=\int e^{-\lambda/\kappa}(\lambda/\kappa)^n/n!\,\mu(d\lambda)\ge0$ and $w_0=r(1/\kappa)$, so the field is the sum of an independent random level with variance $r(1/\kappa)$ and a field whose covariance decays exponentially, rather than a field with slowly decaying covariance. On a general graph, Theorem~\ref{thm:metric_properties}\,(iv) does not imply such a decomposition, and the model differs from a random level added to one short-range covariance.

\subsection{Computation}\label{sec:computation}
Likelihood-based estimation and prediction at $n$ observation locations require the $n\times n$ matrix $D=(d_\kappa(s_i,s_j))_{i,j=1}^{n}$ and the quantities $\log|C|$ and $C^{-1}y$ for $C=r(D)+\sigma_e^2 I$, where $r$ is applied entry-wise and $\sigma_e^2$ is the variance of the measurement noise. Since $C$ is dense, these require $O(n^3)$ operations, which can be reduced by standard approximations such as covariance tapering, with tapers that are positive semidefinite on the graph, low-rank methods or Vecchia approximations. The step that is specific to metric graphs is the computation of $D$, which would require $O(N^3)$ operations and $O(N^2)$ storage if based on a dense inverse or pseudoinverse of a Laplacian matrix on $N$ vertices.

Proposition~\ref{prop:dkappa_computation} reduces this cost through sparsity. On the subdivision $\Gamma(s_1,\ldots,s_n)$, with vertex set $\mathcal V^\star$ of size $O(|\mathcal V|+n)$, we have $d_\kappa(s_i,s_j)=\Sigma_{ii}+\Sigma_{jj}-2\Sigma_{ij}$ with $\Sigma=(Q^\kappa)^{-1}$, so $D$ is obtained from the $n$ columns of $\Sigma$ at the observation locations, that is, from $n$ solves with the sparse matrix $Q^\kappa$. Assembling $Q^\kappa$ costs $O(|\mathcal E|)$ operations, and for planar graphs, or more generally graphs with separators of size $O(|\mathcal V^\star|^{1/2})$, a Cholesky factor of $Q^\kappa$ computed with a nested dissection ordering \citep{lipton_rose_tarjan1979} has $O(|\mathcal V^\star|\log|\mathcal V^\star|)$ non-zero elements and requires $O(|\mathcal V^\star|^{3/2})$ operations. Each solve then costs $O(|\mathcal V^\star|\log|\mathcal V^\star|)$, so that $D$ is computed in $O(|\mathcal V^\star|^{3/2}+n|\mathcal V^\star|\log|\mathcal V^\star|)$ operations. Street and river networks are often close to planar, but these orders rely on the separator assumption. For $\kappa=0$, the matrix $Q$ is singular with kernel spanned by $\mathbf 1$ by Proposition~\ref{prop:Q_properties}\,(ii), and the same computations are done with the positive definite matrix obtained from $Q$ by removing the row and column of one reference vertex, which does not change the resulting distances. 

Prediction at $m$ further locations does not require a new factorisation: since $u_\kappa$ is a Markov process \citep{BSW_Markov}, its value at a point in the interior of an edge depends on the vertex values of $\Gamma(s_1,\ldots,s_n)$ only through the two endpoints of that edge, with conditional mean and variance available in closed form. The distances involving the new locations are therefore obtained from the columns of $\Sigma$ at these endpoints, which requires at most $2m$ further solves with the existing Cholesky factor, together with the closed-form covariance of $u_\kappa$ on an edge given its endpoint values for new locations on the same edge.
The sparsity is also useful for simulation. The algorithms of \citet{alegria2026simulation} for isotropic Gaussian processes on graphs with Euclidean edges require samples of the auxiliary process $Z_\Gamma$ in \eqref{eq:aux_process_resmetric}. 
If the covariance matrix $L^{-1}$ in the definition of $Z_\mu$ is formed explicitly, this step has the usual $O(|\mathcal V^\star|^3)$ dense factorisation cost and $O(|\mathcal V^\star|^2)$ storage. However, 
since the precision matrix $L=Q+e_{v_0}e_{v_0}^\T$ of the vertex values of $Z_\mu$ is sparse, these can be sampled using a sparse Cholesky factor of $L$, at the costs given above, and the Brownian bridges are then simulated edge by edge. 

\section{Application: traffic on two road networks}\label{sec:application}

\subsection{Data and models}

\begin{figure}
    \begin{center}
    \includegraphics[width=0.49\linewidth]{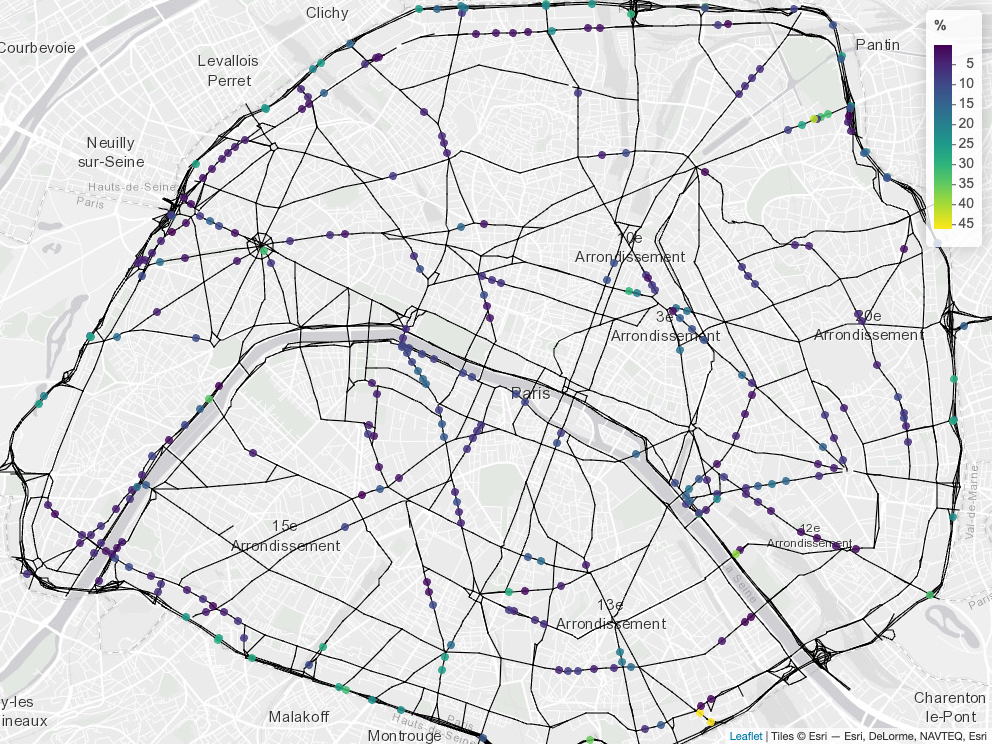}\hfill
    \includegraphics[width=0.49\linewidth]{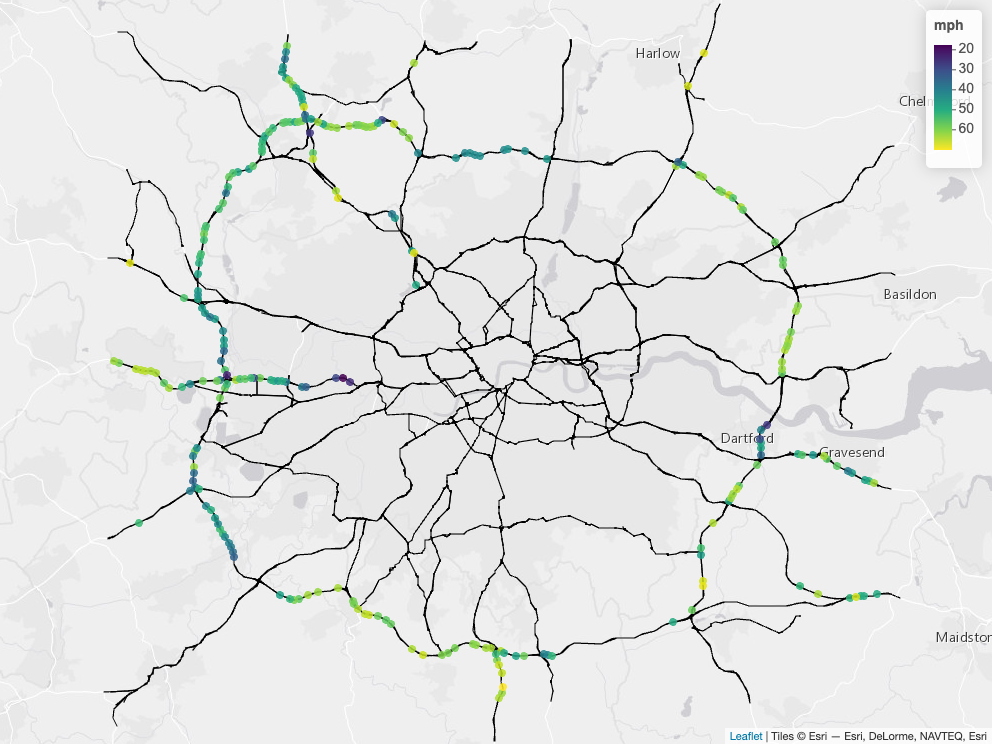}
    \end{center}
    \caption{Road networks and observation locations, coloured by the average over the observed days: occupancy rate in Paris, in \%, averaged on the logarithmic scale (left), and speed around London, in mph (right).}
    \label{fig:data}
\end{figure}

We illustrate the models on two data sets of traffic measurements from inductive loop detectors. The first comes from the permanent traffic counters of the City of Paris \citep{paris_comptages}; we take as observations the logarithm of the occupancy rate, the percentage of time during which a vehicle is above the loop, from $8{:}00$ to $9{:}00$ on working days between 4~September 2025 and 16~October 2026. The network is the largest connected component of the counted road segments, which cover the Boulevard P\'eriph\'erique and the main streets, with one edge for every pair of junctions joined by a segment. Keeping the segments with values on at least $90\%$ of the days and removing days and segments until no values are missing leaves $n=314$ locations, at the segment midpoints, observed on $T=144$ days. With these inserted as vertices, the graph has 2120 vertices, 2640 edges and a total length of 478~km (Figure~\ref{fig:data}). 

The second data set contains speeds from the loop detectors of National Highways on the motorways around London \citep{webtris}. Specifically, the mean speed, in miles per hour, from $17{:}00$ to $18{:}00$ on each day from 1~January to 30~June 2025, at sites at least $0.3$~km apart in the region $0.65^\circ$W--$0.40^\circ$E, $51.22^\circ$N--$51.78^\circ$N. Keeping the sites with values on at least $95\%$ of the days, and the days on which all of them have values, leaves $n=256$ locations on the M25 and the radial motorways, observed on $T=88$ days. The network consists of the motorways and trunk roads in the region with their link roads, from OpenStreetMap \citep{openstreetmap}. The largest connected component, with the locations inserted as vertices and the other vertices of degree two removed, has 3990 vertices, 5207 edges and a total length of 2488~km. 

For both data sets we treat the daily maps as independent replicates of a Gaussian field with common mean and covariance parameters, and compare nine models. All have independent Gaussian measurement errors with variance $\sigma_e^2$ and a mean that is linear in the easting and northing of $s$, and they differ only in the covariance of the latent process. Six are isotropic models with exponential, Mat\'ern with $\nu\le\nicefrac12$, or generalised Cauchy covariance functions \citep[see][Table 1]{anderes2020isotropic}, combined with either the metric $d_0$ or the metric $d_\kappa$ with $\kappa$ estimated from the data; all of these are valid by Theorem~\ref{thm:cm_positive_def}. As comparisons we include the Whittle--Mat\'ern field of \citet{BSW2022} with $\alpha=1$ and a finite element approximation of a Whittle--Mat\'ern field for which $\alpha$ is estimated \citep{Lenin25}, using a mesh with edges of length at most $0.25$~km, a rational approximation of order 2, and $\alpha\in(\nicefrac12,\nicefrac52]$. Finally, as a numerical baseline we include an exponential covariance in the distance $d_{\mathrm{res}}$ obtained by applying \eqref{eq:res_metric_def}--\eqref{eq:aux_process_resmetric} directly to the supplied representation with the locations as vertices, that is, with weights $1/d(v_i,v_j)$ and parallel edges merged, as in \texttt{MetricGraph} \citep{MetricGraph}. As its edges are not Euclidean, $d_{\mathrm{res}}$ differs from $d_0$, by up to $0.006$ and $0.18$~km for the two data sets, and the model is not covered by \citet{anderes2020isotropic}, although its covariance matrices were positive definite in all our fits. All parameters are estimated by maximum likelihood, with the regression coefficients profiled out by generalised least squares. 

Table~\ref{tab:application_results} reports the maximised log-likelihoods. For both data sets, $\kappa$ is estimated to be clearly positive, and $d_\kappa$ increases the maximised log-likelihood relative to $d_0$ for each covariance function. 
For Paris, $\hat\kappa=0.95$~km$^{-1}$ for the exponential and Mat\'ern models, which coincide at $\hat\nu=\nicefrac12$, and $0.59$~km$^{-1}$ for the generalised Cauchy model, a screening length of a few blocks of the street grid. 
For London, $\hat\kappa$ is between $0.14$ and $0.15$~km$^{-1}$, a screening length of about 7~km. 
The generalised Cauchy model in $d_\kappa$ has the highest likelihood for both data sets and the Whittle--Mat\'ern models, whose marginal variances vary along the network, the lowest; refining the finite element mesh to $0.1$~km does not change their rank. 

\subsection{Cross-validation}

To assess prediction at new locations we use leave-percentage-out cross-validation, in which all observations at a held-out location are removed, all parameters, including $\kappa$ and $\alpha$, are re-estimated by maximum likelihood from the training locations, and the held-out values are predicted by their conditional distribution. We hold out $75\%$ or $90\%$ of the locations, drawn at random, 25 times for each percentage with the same draws for all models; the median distance along the network to the nearest training location is then 0.5 and 1.0~km for Paris and 2.5 and 5.1~km for London. We report the root mean squared error (RMSE) and the continuous ranked probability score (CRPS) \citep{gneiting_raftery2007}, averaged over the draws.

\begin{table}
\caption{Number of parameters $p$, including the three mean coefficients, maximised negative log-likelihood, and leave-percentage-out cross-validation: mean RMSE and CRPS over 25 random draws, for the logarithm of the occupancy rate in Paris and in mph for London. Smaller is better; the best value in each column is in bold.}\label{tab:application_results}
\centering
{\small
\begin{tabular*}{\linewidth}{@{\extracolsep{\fill}}lrrrrrr@{}}
\hline
 & & & \multicolumn{2}{c}{75\% held out} & \multicolumn{2}{c}{90\% held out}\\
Model & $p$ & $-\log L$ & RMSE & CRPS & RMSE & CRPS\\
\hline
\multicolumn{7}{@{}l}{\emph{Paris}}\\
$r(d_0)$, exponential & 6 & $62758.307$ & $1.046$ & $0.579$ & $1.083$ & $0.602$\\
$r(d_0)$, generalised Cauchy & 8 & $62593.873$ & $1.041$ & $0.576$ & $1.078$ & $0.600$\\
$r(d_0)$, Mat\'ern & 7 & $62725.681$ & $1.047$ & $0.579$ & $1.079$ & $0.600$\\
\hline
$r(d_\kappa)$, exponential & 7 & $62546.029$ & $1.042$ & $0.577$ & $1.080$ & $0.601$\\
$r(d_\kappa)$, generalised Cauchy & 9 & $\mathbf{62513.985}$ & $\mathbf{1.040}$ & $\mathbf{0.575}$ & $\mathbf{1.076}$ & $\mathbf{0.599}$\\
$r(d_\kappa)$, Mat\'ern & 8 & $62546.029$ & $1.042$ & $0.577$ & $\mathbf{1.076}$ & $\mathbf{0.599}$\\
\hline
Whittle--Mat\'ern $\alpha=1$ & 6 & $63147.079$ & $1.050$ & $0.582$ & $1.092$ & $0.609$\\
Whittle--Mat\'ern FEM, estimated $\alpha$ & 7 & $63059.904$ & $1.051$ & $0.583$ & $1.098$ & $0.614$\\
$r(d_{\mathrm{res}})$, exponential & 6 & $62756.287$ & $1.046$ & $0.579$ & $1.083$ & $0.602$\\
\hline
\multicolumn{7}{@{}l}{\emph{London}}\\
$r(d_0)$, exponential & 6 & $83474.478$ & $11.182$ & $6.076$ & $12.633$ & $6.910$\\
$r(d_0)$, generalised Cauchy & 8 & $83440.611$ & $11.140$ & $\mathbf{6.055}$ & $12.617$ & $6.903$\\
$r(d_0)$, Mat\'ern & 7 & $83448.580$ & $11.176$ & $6.072$ & $12.624$ & $6.905$\\
\hline
$r(d_\kappa)$, exponential & 7 & $83321.751$ & $11.141$ & $6.057$ & $\mathbf{12.581}$ & $\mathbf{6.886}$\\
$r(d_\kappa)$, generalised Cauchy & 9 & $\mathbf{83317.492}$ & $\mathbf{11.139}$ & $\mathbf{6.055}$ & $\mathbf{12.581}$ & $\mathbf{6.886}$\\
$r(d_\kappa)$, Mat\'ern & 8 & $83317.907$ & $11.140$ & $\mathbf{6.055}$ & $\mathbf{12.581}$ & $\mathbf{6.886}$\\
\hline
Whittle--Mat\'ern $\alpha=1$ & 6 & $83543.142$ & $\mathbf{11.139}$ & $6.085$ & $12.638$ & $6.921$\\
Whittle--Mat\'ern FEM, estimated $\alpha$ & 7 & $83869.791$ & $11.183$ & $6.132$ & $12.613$ & $6.918$\\
$r(d_{\mathrm{res}})$, exponential & 6 & $83475.250$ & $11.182$ & $6.076$ & $12.633$ & $6.910$\\
\hline
\end{tabular*}}
\end{table}

In Table~\ref{tab:application_results}, the $d_\kappa$ models have the smallest RMSE and CRPS for both data sets and both percentages, with ties for London with $75\%$ held out. For London, $d_\kappa$ reduces the RMSE of the exponential model by $0.04$ and $0.05$~mph and its CRPS by $0.3\%$ and $0.4\%$, with the smaller CRPS in 19 and 18 of the 25 draws. For Paris, the generalised Cauchy model in $d_\kappa$ has scores $0.5\%$ to $0.6\%$ below those of the exponential model in $d_0$, and the exponential model in $d_\kappa$ has a smaller CRPS than that in $d_0$ in 15 and 17 of the 25 draws. The median estimates of $\kappa$ in the training sets are $0.65$ and $0.50$~km$^{-1}$ for Paris and $0.10$ and $0.05$~km$^{-1}$ for London.

Thus, for both data sets, $\kappa$ is estimated to be clearly positive, $d_\kappa$ improves the likelihood for every covariance function, and it gives the best predictions at locations held out at random. The differences in prediction are below $1\%$, as the predictions are dominated by the nearest training locations, at distances for which $d_\kappa$ and $d_0$ are similar.

\section{Discussion}\label{sec:discussion_top}


While this manuscript was being finalised, \citet{filosi2026polyharmonic} independently introduced  polyharmonic distances on metric graphs, whose squares are the limits as $\kappa\downarrow0$ of the increment variances of our Whittle--Mat\'ern fields of integer order $p$. For $p=1$ this is $d_0$, and their results on its invariance under subdivision, its agreement with the resistance metric of \citet{anderes2020isotropic} and the validity of completely monotone functions of it overlap with ours for $\kappa=0$; they are formulated for simple graphs, with other representations handled by inserting vertices. They do not consider $\kappa>0$, which is the source of the long-range dependence studied here, or the computation of the distances. A natural extension is to combine the two constructions through the increment variances of higher-order Whittle--Mat\'ern fields with $\kappa>0$. Composing these with completely monotone functions gives positive semidefinite kernels, but whether they have a sparse finite-dimensional representation like that of the first-order case requires a separate analysis.

Beyond Gaussian processes, the matrix $Q$ may be useful as a length-weighted Laplacian for graph signal processing and graph learning on networks, where consistency under subdivision is a natural requirement for methods that should not depend on the representation for a given network, although methods built on $Q$ do not automatically inherit it.

\section*{Acknowledgement}
This research was funded by King Abdullah University of Science and Technology (KAUST) under Award No. ORFS-CRG12-2024-6399.

\appendix
\setcounter{equation}{0}
\renewcommand{\theequation}{A\arabic{equation}}
\section{Proofs}\label{sec:proofs}

\subsection{Preliminaries}\label{sec:proofs_preliminaries}

The following proposition gives the one-dimensional energy-minimisation identity used repeatedly in the proofs. It connects the continuum energy of $\kappa^2-\Delta_\Gamma$ with the vertex matrices $Q^\kappa$ and $Q$, and is used in Proposition~\ref{prop:Q_properties}(vii) and Theorem~\ref{thm:metric_properties}(i)--(ii).

\begin{proposition}\label{prop:edge_DN}
Let $I=[0,\ell]$, with $\ell>0$, and for $\kappa\ge0$ define
\[
\mathcal E_\kappa(f)=\int_0^\ell\{|f'(x)|^2+\kappa^2|f(x)|^2\}\,dx,
\qquad f\in H^1(I).
\]
Fix endpoint values $a,b\in\mathbb R$ and minimise $\mathcal E_\kappa(f)$ over $f\in H^1(I)$ satisfying $f(0)=a$ and $f(\ell)=b$.
\begin{enumerate}
\item[(i)] If $\kappa=0$, the unique minimiser is
$$
h_0(x)=a+\frac{b-a}{\ell}x,
\quad
\text{and}
\quad
\min \mathcal E_0=\frac{(a-b)^2}{\ell}.
$$
\item[(ii)] If $\kappa>0$, the unique minimiser is
\begin{equation}\label{eq:edge_DN_appendix}
h_\kappa(x)=a\frac{\sinh\{\kappa(\ell-x)\}}{\sinh(\kappa\ell)}
+b\frac{\sinh(\kappa x)}{\sinh(\kappa\ell)},
\quad
\text{and}
\quad
\min \mathcal E_\kappa
=\kappa\coth(\kappa\ell)(a^2+b^2)
 -2\kappa\operatorname{csch}(\kappa\ell)ab.
\end{equation}
\item[(iii)] Consequently, on a finite metric graph with prescribed vertex-value vector $x$, minimising the sum of the edge energies over continuous edge functions with those vertex traces gives $x^TQx$ when $\kappa=0$ and $x^TQ^\kappa x$ when $\kappa>0$.  In particular, a self-loop based at a vertex with value $a$ contributes zero when $\kappa=0$ and $2\kappa\tanh(\kappa\ell/2)a^2$ when $\kappa>0$.
\end{enumerate}
\end{proposition}

\noindent\emph{Proof.}
For $\kappa=0$, write $f=h_0+g$ with $g\in H^1_0(I)$. Since $h_0'=(b-a)/\ell$ is constant,
$\int_0^\ell h_0'(x)g'(x)\,dx
=\frac{b-a}{\ell}\{g(\ell)-g(0)\}=0$.
Therefore
$\mathcal E_0(f)=\frac{(a-b)^2}{\ell}+\int_0^\ell|g'(x)|^2\,dx$,
which proves (i), including uniqueness.

For $\kappa>0$, the displayed function $h_\kappa$ is the unique solution of
$-h_\kappa''+\kappa^2h_\kappa=0$ with the prescribed endpoint values.  If $f=h_\kappa+g$ with $g\in H^1_0(I)$, then integration by parts gives
$\int_0^\ell(h_\kappa'g'+\kappa^2h_\kappa g)\,dx
=[h_\kappa'g]_0^\ell
 +\int_0^\ell(-h_\kappa''+\kappa^2h_\kappa)g\,dx=0$.
Hence
$\mathcal E_\kappa(f)=\mathcal E_\kappa(h_\kappa)+\mathcal E_\kappa(g)$,
so $h_\kappa$ is the unique minimiser.  Using
$\mathcal E_\kappa(h_\kappa)=[h_\kappa h_\kappa']_0^\ell$ and the endpoint derivatives of $h_\kappa$ yields \eqref{eq:edge_DN_appendix}.

For (iii), once the vertex values are fixed, the edge interiors are independent in the energy minimisation, so the total minimum is the sum of the edgewise minima from (i)--(ii). For a non-loop edge with endpoint values $x_i,x_j$, these are respectively
$\ell_e^{-1}(x_i-x_j)^2$
and
$\kappa\coth(\kappa\ell_e)(x_i^2+x_j^2)
 -2\kappa\operatorname{csch}(\kappa\ell_e)x_ix_j$,
which are exactly the corresponding contributions to $x^TQx$ and $x^TQ^\kappa x$. For a self-loop the two endpoint traces coincide, so setting $a=b$ in (i)--(ii) gives the stated loop contributions, using $\coth z-\operatorname{csch}z=\tanh(z/2)$. Summing over all edges proves (iii).\hfill\ensuremath{\square}

\subsection{Proof of Proposition~\ref{prop:Q_properties}}

Symmetry is immediate from \eqref{eq:Q0}. Expanding $x^\T Q x$ and grouping the contribution of each non-loop edge $e=(v_i,v_j)$ gives
$\frac{x_i^2+x_j^2-2x_ix_j}{\ell_e}=\frac{(x_i-x_j)^2}{\ell_e}$,
while a loop contributes zero. This proves (i), and (ii) follows because every term is nonnegative. If $x^\T Q x=0$, then $x_i=x_j$ across every non-loop edge. Connectivity of the underlying metric graph implies that all vertex values agree, so $\ker Q=\operatorname{span}\{\mathbf1\}$.

Part (iii) follows by comparing \eqref{eq:Q0} with \eqref{eq:standard_laplacian}. Part (iv) follows because every parallel edge contributes $-1/\ell_e$ to the same off-diagonal entry , and the corresponding conductances add. Part (v) follows directly from the indicator in the diagonal and from the fact that a loop has no pair of distinct endpoints.

For (vi), if $e=(v_i,v_j)$ is non-loop, its contribution to $Q^\kappa$ is
\[
\kappa\begin{pmatrix}
\coth(\kappa\ell_e)&-\operatorname{csch}(\kappa\ell_e)\\
-\operatorname{csch}(\kappa\ell_e)&\coth(\kappa\ell_e)
\end{pmatrix}.
\]
Using $z\coth z\to1$ and $z\operatorname{csch}z\to1$ as $z\downarrow0$, this converges to
\mbox{$\ell_e^{-1}\left(\begin{smallmatrix}1&-1\\-1&1\end{smallmatrix}\right)$
}. A loop contributes $2\kappa\tanh(\kappa\ell_e/2)\to0$. Summing over finitely many edges proves entrywise convergence.

For (vii), let $x\in\mathbb R^{|\mathcal V|}$ be prescribed on the original vertices and let $z\in\mathbb R^{|\mathcal I|}$ denote values on inserted vertices. By Proposition~\ref{prop:edge_DN}(i),(iii),
\[
\begin{pmatrix}x\\z\end{pmatrix}^{\!\T}\widetilde Q
\begin{pmatrix}x\\z\end{pmatrix}
\]
is the minimum Dirichlet energy among continuous edgewise-affine functions on the refined graph with those nodal values. Minimising first over $z$ is equivalent to removing the artificial subdivision constraints and therefore gives the minimum energy associated with the original nodal vector $x$, namely $x^\T Q x$.

Write
\mbox{
$\widetilde Q=\left(\begin{smallmatrix}A&B\\B^\T&C\end{smallmatrix}\right)$
}
according to $\mathcal V\cup\mathcal I$. If $z^TCz=0$, extend $z$ by zero on $\mathcal V$. The quadratic-form identity in (i) then forces this extended vector to be constant on the connected refined graph, hence zero. Thus $C\succ0$. The minimiser in $z$ is hence $z_*=-C^{-1}B^Tx$, and substitution gives
$x^TQx=x^T(A-BC^{-1}B^T)x$ 
for every $x$. Equality of the symmetric matrices proves \eqref{eq:kron_subdivision}. \hfill\ensuremath{\square}

\subsection{Proof of Proposition~\ref{prop:dkappa_computation}}

Fix the subdivision $\Gamma(s_1,s_2)$ and write $e=e_{s_1s_2}$. For $\kappa>0$, the vertex vector of $u_\kappa$ is a centred Gaussian vector with precision $Q^\kappa$, hence
$d_\kappa(s_1,s_2)=\var(e^\T u_\kappa)=e^\T(Q^\kappa)^{-1}e$.
The matrix $Q^\kappa$ is positive definite: using that 
$\coth z=\operatorname{csch}z+\tanh(z/2)$,
each non-loop edge contributes
$c_e(x_i-x_j)^2+g_e(x_i^2+x_j^2)$, with 
$c_e=\kappa\operatorname{csch}(\kappa\ell_e)$ and 
$g_e=\kappa\tanh(\kappa\ell_e/2)$,
and a loop contributes $2g_ex_i^2$. Since $c_e,g_e>0$ and every vertex is incident to an edge, the total quadratic form is strictly positive for $x\ne0$.

It remains to identify the $\kappa\downarrow0$ limit. For $\kappa>0$, set
\[
q_\kappa(f,g)=\sum_{e\in\mathcal E}\int_e
\{f'(x)g'(x)+\kappa^2f(x)g(x)\}\,dx,
\qquad f,g\in H^1(\Gamma).
\]
This is the closed coercive quadratic form associated with the operator $\kappa^2-\Delta_\Gamma$. Since point evaluation is continuous on $H^1(\Gamma)$, the contrast functional $\ell_{s_1,s_2}(f)=f(s_1)-f(s_2)$ is continuous with respect to the $q_\kappa$-norm. Hence the Riesz representation theorem gives a unique $g_\kappa\in H^1(\Gamma)$ such that
$q_\kappa(g_\kappa,f)=\ell_{s_1,s_2}(f),
\qquad f\in H^1(\Gamma)$.
The covariance operator of $u_\kappa$ is $(\kappa^2-\Delta_\Gamma)^{-1}$, so $g_\kappa$ is precisely the covariance representer of this contrast and
$d_\kappa(s_1,s_2)=\ell_{s_1,s_2}(g_\kappa)=q_\kappa(g_\kappa,g_\kappa)$. Therefore, by Cauchy--Schwarz in the $q_\kappa$ inner product,
\[
|\ell_{s_1,s_2}(f)|^2
=|q_\kappa(g_\kappa,f)|^2
\le d_\kappa(s_1,s_2)q_\kappa(f,f),
\]
with equality for $f=g_\kappa$. This proves the variational identity
\begin{equation}\label{eq:dual_variational_dkappa}
d_\kappa(s_1,s_2)
=\sup_{f\in H^1(\Gamma),\,f\ne0}
\frac{|f(s_1)-f(s_2)|^2}
{\int_\Gamma\{|f'(x)|^2+\kappa^2|f(x)|^2\}\,dx}.
\end{equation}
As $\kappa\downarrow0$, the denominators decrease pointwise to the Dirichlet energy. Hence the right-hand side increases to
\[
R_0(s_1,s_2)=
\sup_{f:\,\int|f'|^2>0}
\frac{|f(s_1)-f(s_2)|^2}{\int_\Gamma|f'(x)|^2\,dx}.
\]
This supremum is finite: if $\gamma$ is any path from $s_1$ to $s_2$, then Cauchy--Schwarz gives
$|f(s_1)-f(s_2)|^2\le \operatorname{length}(\gamma)\int_\gamma|f'|^2\le \operatorname{length}(\gamma)\mathcal E(f)$, and hence $R_0(s_1,s_2)\le d(s_1,s_2)$.
Indeed, the limit is at most this supremum, while the reverse inequality follows by fixing an $f$ whose zero-mass quotient is within $\varepsilon$ of the supremum and then letting $\kappa\downarrow0$. Thus $d_0=R_0$.

Insert $s_1,s_2$ as vertices. Since the numerator in the variational quotient depends only on the vertex values, Proposition~\ref{prop:edge_DN}(i),(iii) reduces the zero-mass variational problem to
\[
R_0(s_1,s_2)
=
\sup_{x:\,x^TQx>0}
\frac{(e^Tx)^2}{x^TQx},
\qquad e=e_{s_1s_2}.
\]
Because $e^T\mathbf1=0$ and $\ker Q=\operatorname{span}\{\mathbf1\}$, we have $e\in\operatorname{range}Q$. Put $y=Q^\dagger e$, so that $Qy=e$. Then, by Cauchy--Schwarz for the positive-semidefinite bilinear form $(a,b)\mapsto a^TQb$ on $\mathbf1^\perp$,
\[
(e^Tx)^2=(y^TQx)^2
\le (y^TQy)(x^TQx)
=(e^TQ^\dagger e)(x^TQx).
\]
Equality is attained for $x=y$, and hence
$R_0(s_1,s_2)=e^TQ^\dagger e$.
The equation $Qz=e$ is therefore solvable, every solution is of the form $z=Q^\dagger e+c\mathbf1$, and $e^Tz=e^TQ^\dagger e$. For $\kappa>0$, the same finite-dimensional variational argument with the positive-definite matrix $Q^\kappa$ gives $d_\kappa(s_1,s_2)=e^T(Q^\kappa)^{-1}e$.\hfill\ensuremath{\square}

\subsection{Proof of Theorem~\ref{thm:metric_properties}}

\noindent(ii)
Let $A=-\Delta_\Gamma$ with Kirchhoff conditions. Since $\Gamma$ is compact and connected, $A$ is non-negative and self-adjoint with compact resolvent. Write $A\phi_j=\lambda_j\phi_j$, where $0=\lambda_0<\lambda_1\le\lambda_2\le\cdots$, the eigenfunctions form an orthonormal basis of $L^2(\Gamma)$, and $\phi_0$ is constant. Therefore, for $\kappa>0$, the covariance operator $(A+\kappa^2)^{-1}$ gives
\begin{equation}\label{eq:spectral_dkappa}
d_\kappa(s,t)=\sum_{j=1}^{\infty}
\frac{|\phi_j(s)-\phi_j(t)|^2}{\lambda_j+\kappa^2}.
\end{equation}
Proposition~\ref{prop:dkappa_computation} already shows that $d_\kappa(s,t)\to d_0(s,t)<\infty$ as $\kappa\downarrow0$. Since every summand in \eqref{eq:spectral_dkappa} increases to its value with denominator $\lambda_j$, monotone convergence gives
$d_0(s,t)=\sum_{j=1}^{\infty}\frac{|\phi_j(s)-\phi_j(t)|^2}{\lambda_j}$.
The series on the right is finite and dominates the corresponding series for every $\kappa\ge0$. Dominated convergence therefore gives continuity in $\kappa\in[0,\infty)$ and $d_\kappa(s,t)\to0$ as $\kappa\to\infty$. If $s\ne t$, the functional $f\mapsto f(s)-f(t)$ is non-zero and continuous on $H^1(\Gamma)$. Finite linear combinations of the eigenfunctions are dense in the form domain $H^1(\Gamma)$; hence, if $\phi_j(s)=\phi_j(t)$ for every $j\ge1$, the contrast would vanish on a dense subset and therefore on all of $H^1(\Gamma)$, a contradiction. Thus at least one coefficient $\phi_j(s)-\phi_j(t)$ with $j\ge1$ is non-zero. Every term in \eqref{eq:spectral_dkappa} is non-increasing in $\kappa$, and at least one is strictly decreasing. 

\noindent(i)
Symmetry and $d_\kappa(s,s)=0$ are immediate. To prove the triangle inequality, fix three points $s,t,r\in\Gamma$ and work on the single subdivision $\Gamma(s,t,r)$. For $\kappa=0$, Proposition~\ref{prop:dkappa_computation} identifies all three quantities $d_0(s,t)$, $d_0(s,r)$ and $d_0(r,t)$ with resistance distances in the same connected weighted graph, whose edge weights are reciprocal segment lengths. The standard resistance metric is strictly positive off the diagonal and satisfies the triangle inequality \citep{klein_randic1993}.

For $\kappa>0$, the same conclusion can be reduced to the preceding graph-theoretic result without using any physical interpretation, again on the common subdivision $\Gamma(s,t,r)$. For every non-loop edge $e$, set
$c_e=\kappa\operatorname{csch}(\kappa\ell_e)$,
and 
 $g_e=\kappa\tanh(\kappa\ell_e/2)$.
Adjoin one auxiliary vertex $\partial$. Give the original endpoints of $e$ an edge of weight $c_e$, and connect each endpoint to $\partial$ by an edge of weight $g_e$; for a self-loop based at $v$, use one edge of weight $2g_e$ between $v$ and $\partial$. Let $\widehat Q^\kappa$ be the weighted graph Laplacian of this augmented graph, and let $h$ be the vector whose $i$th component is the total weight of the edges joining $v_i$ to $\partial$. Then
\[
\widehat Q^\kappa
=
\begin{pmatrix}
Q^\kappa & -h\\
-h^T & \mathbf1^Th
\end{pmatrix},
\qquad
Q^\kappa\mathbf1=h.
\]
Fix $s,t$ and write $b=e_s-e_t$ and $x=(Q^\kappa)^{-1}b$. Since $\mathbf1^Tb=0$,
$h^Tx=\mathbf1^TQ^\kappa x=\mathbf1^Tb=0$,
and therefore
\[
\widehat Q^\kappa
\begin{pmatrix}x\\0\end{pmatrix}
=
\begin{pmatrix}b\\0\end{pmatrix}.
\]
The vector $(b,0)$ is orthogonal to the constants, so $(x,0)$ differs from $\widehat Q^{\kappa\dagger}(b,0)$ only by an additive constant vector. Consequently,
$(b,0)^T\widehat Q^{\kappa\dagger}(b,0)
=b^Tx
=b^T(Q^\kappa)^{-1}b
=d_\kappa(s,t)$.
Thus $d_\kappa$ is the resistance distance of a connected weighted graph, and is therefore a metric.

For $\kappa>0$, subdivision invariance follows directly from the variational characterization. If an edge $[0,\ell]$ is split at $a\in(0,\ell)$, restriction gives the canonical isometric identification
$H^1(0,\ell)\cong
\{(f_1,f_2)\in H^1(0,a)\oplus H^1(a,\ell):f_1(a)=f_2(a)\}$.
Under this identification all point evaluations are preserved and
\[
\int_0^\ell (|f'|^2+\kappa^2|f|^2)
=
\int_0^a (|f_1'|^2+\kappa^2|f_1|^2)
+
\int_a^\ell (|f_2'|^2+\kappa^2|f_2|^2).
\]
Thus the variational quotient in \eqref{eq:dual_variational_dkappa} defining $d_\kappa$ is unchanged, including when the inserted point itself is one of the evaluation locations. Repeating this argument proves invariance under any finite degree-two subdivision. Since $d_0$ is the pointwise $\kappa\downarrow0$ limit of the subdivision-invariant metrics $d_\kappa$, subdivision invariance also holds at $\kappa=0$.

\noindent(iii)
Proposition~\ref{prop:dkappa_computation} gives the variational characterization
\[
d_0(s,t)
=
\sup_{f:\,\mathcal E(f)>0}
\frac{|f(s)-f(t)|^2}{\mathcal E(f)},
\qquad
\mathcal E(f)=\sum_{e\in\mathcal E}\int_e|f'(x)|^2\,dx.
\]
Let $\gamma$ be a geodesic from $s$ to $t$, parametrised by arc length. For every $f\in H^1(\Gamma)$,
\[
|f(s)-f(t)|^2
=
\left|\int_\gamma f'(x)\,dx\right|^2
\le d(s,t)\int_\gamma|f'(x)|^2\,dx
\le d(s,t)\mathcal E(f),
\]
by Cauchy--Schwarz. Taking the supremum gives $d_0(s,t)\le d(s,t)$.

Suppose first that $\Gamma$ is a metric tree. The path $\gamma$ joining $s$ and $t$ is then unique. Put $L=d(s,t)$ and define $f$ to increase linearly with slope one along $\gamma$, with $f(s)=0$, and to be constant on every component attached to $\gamma$. This gives a continuous function in $H^1(\Gamma)$ satisfying
$f(t)-f(s)=L$,
$\mathcal E(f)=L$.
Hence $d_0(s,t)\ge L=d(s,t)$, and equality holds for every pair.

Conversely, suppose that $\Gamma$ contains a cycle. Choose an edge belonging to that cycle and two distinct interior points $u,v$ on the edge sufficiently close that the direct subsegment between them, denoted by $\gamma_1$, is the unique geodesic. Let its length be $\ell_1=d(u,v)$, and let $\gamma_2$ denote the other arc of the cycle, with length $\ell_2>0$. The interiors of $\gamma_1$ and $\gamma_2$ are disjoint. For every $f\in H^1(\Gamma)$, Cauchy--Schwarz on the two paths gives
\[
\int_{\gamma_i}|f'(x)|^2\,dx
\ge
\frac{|f(u)-f(v)|^2}{\ell_i},
\qquad i=1,2.
\]
Therefore
\[
\mathcal E(f)
\ge
|f(u)-f(v)|^2
\left(\frac1{\ell_1}+\frac1{\ell_2}\right), \quad 
\text{and hence} \quad
d_0(u,v)
\le
\left(\frac1{\ell_1}+\frac1{\ell_2}\right)^{-1}
<\ell_1=d(u,v).
\]
Thus equality $d_0(s,t)=d(s,t)$ for every pair $s,t$ is possible only when $\Gamma$ is a metric tree. 

\noindent(iv) 
If $s=t$, \eqref{eq:screening_bound} is immediate. Assume $s\ne t$, fix $\kappa>0$, and identify a geodesic $\gamma$ from $s$ to $t$ with the interval $[0,L]$, where $L=d(s,t)>0$. For any $f\in H^1(\Gamma)$, let $a=f(s)$ and $b=f(t)$. Since the full energy dominates the energy along $\gamma$, Proposition~\ref{prop:edge_DN}(ii) gives
\[
\int_\Gamma (|f'|^2+\kappa^2|f|^2)
\ge
\kappa\coth(\kappa L)(a^2+b^2)
-2\kappa\operatorname{csch}(\kappa L)ab.
\]
The vector $(1,-1)^T$ is an eigenvector of the associated $2\times2$ matrix with eigenvalue
$\kappa\{\coth(\kappa L)+\operatorname{csch}(\kappa L)\}=\kappa\coth(\kappa L/2)$. Hence
\[
\frac{|a-b|^2}{\kappa\coth(\kappa L)(a^2+b^2)-2\kappa\operatorname{csch}(\kappa L)ab}
\le \frac{2}{\kappa}\tanh\left(\frac{\kappa L}{2}\right).
\]
Taking the supremum in \eqref{eq:dual_variational_dkappa} proves \eqref{eq:screening_bound}; the strict upper bound $2/\kappa$ follows from $\tanh x<1$ for finite $x$.
\hfill\ensuremath{\square}

\subsection{Increasing-domain calculation on the line}

Let $\Gamma_T=[-T,T]$ and let $u_{\kappa,T}$ be the first-order Whittle--Mat\'ern field on this interval. The interval is a single edge whose endpoints are vertices of degree one, so the Kirchhoff conditions reduce to Neumann conditions, and by \citet{BSW2023_AOS} the covariance function of $u_{\kappa,T}$ is
\begin{equation}\label{eq:interval_covariance}
\varrho_{\kappa,T}(s,t)
=\frac{\cosh\{\kappa(2T-|s-t|)\}+\cosh\{\kappa(s+t)\}}
{2\kappa\sinh(2\kappa T)},
\qquad s,t\in[-T,T].
\end{equation}
Since $d_{\kappa,T}$ is the variogram of $u_{\kappa,T}$,
\begin{equation}\label{eq:dkappa_interval_cov}
d_{\kappa,T}(s,t)
=\varrho_{\kappa,T}(s,s)+\varrho_{\kappa,T}(t,t)-2\varrho_{\kappa,T}(s,t).
\end{equation}
For fixed $s,t\in\mathbb R$, the second hyperbolic-cosine term in the numerator of \eqref{eq:interval_covariance} is independent of $T$, whereas
$\frac{\cosh\{\kappa(2T-|s-t|)\}}{2\kappa\sinh(2\kappa T)}
\longrightarrow \frac{1}{2\kappa}e^{-\kappa|s-t|}$.
Consequently,
$\varrho_{\kappa,T}(s,t)\longrightarrow \frac{1}{2\kappa}e^{-\kappa|s-t|}$,
which is the covariance function of the stationary Ornstein--Uhlenbeck process on $\mathbb R$, equivalently the first-order Whittle--Mat\'ern field on the line under the normalization used here. Substitution into \eqref{eq:dkappa_interval_cov} yields
$\lim_{T\to\infty}d_{\kappa,T}(s,t)
=\frac1\kappa\{1-e^{-\kappa|s-t|}\}$ for $\kappa>0$.
For $\kappa=0$, the interval is a metric tree, and Theorem~\ref{thm:metric_properties}(iii) gives $d_{0,T}(s,t)=|s-t|$ for every $T>\max(|s|,|t|)$. This is the variogram of a two-sided standard Brownian motion on $\mathbb R$, since $\operatorname{var}\{B(s)-B(t)\}=|s-t|$. Combining the two cases proves \eqref{eq:line_metric_limit}.

\subsection{Proof of Theorem~\ref{thm:resistance_equiv}}

Assume that $\Gamma$ has Euclidean edges. Then it is simple and loop-free, and every edge $e=(v_i,v_j)$ satisfies $d(v_i,v_j)=\ell_e$. Hence the conductance in the resistance construction is $1/d(v_i,v_j)=1/\ell_e$, so its singular conductance Laplacian is exactly $Q$.

For arbitrary $s,t\in\Gamma$, insert them as degree-two vertices. The refined graph still has Euclidean edges: if a refined segment of an original edge admitted a shorter external path between its endpoints, adjoining the complementary pieces of the original edge would produce a path between the original endpoints shorter than the original edge, contradicting the Euclidean-edge assumption. On this refined graph, the resistance metric of \citet{anderes2020isotropic} is the effective resistance of the conductance network with segment conductance equal to reciprocal segment length. Proposition~\ref{prop:dkappa_computation} identifies $d_0(s,t)$, defined by the $\kappa\downarrow0$ limit, with the effective resistance of exactly the same network. Thus $d_0=d_R$.\hfill\ensuremath{\square}

\subsection{Proof of Theorem~\ref{thm:cm_positive_def}}

Fix $s_1,\ldots,s_m\in\Gamma$ and real coefficients $a_1,\ldots,a_m$ satisfying $\sum_i a_i=0$. Work on the subdivision $\Gamma(s_1,\ldots,s_m)$, in which all the points are vertices. Let $G_\kappa= (Q^\kappa)^{-1}$ if $\kappa>0$ and $G_0 = Q^\dagger$.
Let $G$ denote the principal submatrix of $G_\kappa$ corresponding to $s_1,\ldots,s_m$. Proposition~\ref{prop:dkappa_computation} gives
$d_\kappa(s_i,s_j)=G_{ii}+G_{jj}-2G_{ij}$. Therefore
$\sum_{i,j=1}^m a_ia_jd_\kappa(s_i,s_j)
=-2a^TGa\le0$,
because the diagonal terms vanish when $\sum_i a_i=0$ and $G$ is positive semidefinite. Thus $d_\kappa$ is conditionally negative definite.

By Bernstein's theorem, a completely monotone function has a Laplace-transform representation
$r(t)=\int_{[0,\infty)}e^{-\lambda t}\,\mu(d\lambda)$ for a finite positive measure $\mu$. Schoenberg's theorem implies that $e^{-\lambda d_\kappa(s,t)}$ is positive semidefinite for every $\lambda\ge0$ whenever $d_\kappa$ is conditionally negative definite and vanishes on the diagonal. Integrating these positive-semidefinite kernels against $\mu$ proves that $r(d_\kappa(s,t))$ is positive semidefinite. \hfill\ensuremath{\square}

\bibliographystyle{plainnat}
  \bibliography{references_final_rigorous}

@unpublished{Lenin25,
  author = {David Bolin and Lenin Riera-Segura and Alexandre B. Simas},
  note = {preprint 	arXiv:2501.11738},
  title = {A new class of non-stationary Gaussian fields with general smoothness on metric graphs},
  year = {2025}
}

@article{ROAD2,
	author = {Ramirez, Andres Felipe and Valencia, Carlos},
	journal = {Accid. Anal. Prev.},
	month = {JAN},
	title = {{Spatiotemporal correlation study of traffic accidents with fatalities and injuries in Bogota (Colombia)}},
	volume = {149},
	year = {2021}}

@article{ROAD3,
	author = {Chaudhuri, Somnath and Juan, Pablo and Mateu, Jorge},
	journal = {J. Appl. Statist.},
	month = {DEC 10},
	number = {16},
	pages = {3229-3250},
	title = {Spatio-temporal modeling of traffic accidents incidence on urban road networks based on an explicit network triangulation},
	volume = {50},
	year = {2023}}

@Manual{MetricGraph,
  title = {{\textit{MetricGraph}}: Random fields on metric graphs},
  author = {David Bolin and Alexandre B. Simas and Jonas Wallin},
  year = {2023},
  note = {R package version 1.4.0},
  url = {https://CRAN.R-project.org/package=MetricGraph},
}

@incollection{baker_faber2006,
  author    = {Baker, Matthew and Faber, Xander},
  title     = {Metrized graphs, {L}aplacian operators, and electrical networks},
  booktitle = {Quantum Graphs and Their Applications},
  series    = {Contemporary Mathematics},
  volume    = {415},
  pages     = {15--33},
  publisher = {American Mathematical Society},
  year      = {2006},
}

@article{baker_rumely2007,
  author    = {Baker, Matthew and Rumely, Robert},
  title     = {Harmonic analysis on metrized graphs},
  journal   = {Can. J. Math.},
  volume    = {59},
  number    = {2},
  pages     = {225--275},
  year      = {2007},
}

@article{cinkir2011,
  author  = {Cinkir, Zubeyir},
  title   = {The tau constant and the discrete {L}aplacian matrix of a metrized graph},
  journal = {Eur. J. Combin.},
  volume  = {32},
  number  = {4},
  pages   = {639--655},
  year    = {2011},
}

@book{berkolaiko_kuchment2013,
  author    = {Berkolaiko, Gregory and Kuchment, Peter},
  title     = {Introduction to Quantum Graphs},
  series    = {Mathematical Surveys and Monographs},
  volume    = {186},
  publisher = {American Mathematical Society},
  address   = {Providence, RI},
  year      = {2013},
}

@article{dorfler_bullo2013,
  author  = {D{\"o}rfler, Florian and Bullo, Francesco},
  title   = {{K}ron reduction of graphs with applications to electrical networks},
  journal = {IEEE Trans. Circuits Syst. I},
  volume  = {60},
  number  = {1},
  pages   = {150--163},
  year    = {2013},
}

@article{peterson2013,
  author  = {Peterson, Erin E. and Ver Hoef, Jay M. and Isaak, Daniel J. and Falke, Jeffrey A. and Fortin, Marie-Jos\'ee and Jordan, Chris E. and McNyset, Kristina and Monestiez, Pascal and Ruesch, Aaron S. and Sengupta, Aritra and Som, Nicholas and Steel, E. Ashley and Theobald, David M. and Torgersen, Christian E. and Wenger, Seth J.},
  title   = {Modelling dendritic ecological networks in space: an integrated network perspective},
  journal = {Ecol. Lett.},
  volume  = {16},
  number  = {5},
  pages   = {707--719},
  year    = {2013},
}

@article{baddeley2021review,
  author  = {Baddeley, Adrian and Nair, Gopalan and Rakshit, Suman and McSwiggan, Greg and Davies, Tilman M.},
  title   = {Analysing point patterns on networks --- {A} review},
  journal = {Spatial Statist.},
  volume  = {42},
  pages   = {100435},
  year    = {2021},
}

@inproceedings{kondor_lafferty2002,
  author    = {Kondor, Risi I. and Lafferty, John},
  title     = {Diffusion kernels on graphs and other discrete input spaces},
  booktitle = {Proc. 19th Int. Conf. Mach. Learn.},
  pages     = {315--322},
  publisher = {Morgan Kaufmann},
  year      = {2002},
}

@inproceedings{smola_kondor2003,
  author    = {Smola, Alexander J. and Kondor, Risi},
  title     = {Kernels and regularization on graphs},
  booktitle = {Learning Theory and Kernel Machines: Proc. 16th Annu. Conf. Learn. Theory and 7th Kernel Workshop},
  pages     = {144--158},
  publisher = {Springer},
  year      = {2003},
}

@book{doyle_snell1984,
  author    = {Doyle, Peter G. and Snell, J. Laurie},
  title     = {Random walks and electric networks},
  series    = {Carus Mathematical Monographs},
  volume    = {22},
  publisher = {Mathematical Association of America},
  address   = {Washington, DC},
  year      = {1984},
}

@article{klein_randic1993,
  author    = {Klein, D. J. and Randi{\'c}, M.},
  title     = {Resistance distance},
  journal   = {J. Math. Chem.},
  volume    = {12},
  number    = {1},
  pages     = {81--95},
  year      = {1993},
}

@article{lipton_rose_tarjan1979,
  author    = {Lipton, Richard J. and Rose, Donald J. and Tarjan, Robert Endre},
  title     = {Generalized nested dissection},
  journal   = {SIAM J. Numer. Anal.},
  volume    = {16},
  number    = {2},
  pages     = {346--358},
  year      = {1979},
}

@article{alegria2026simulation,
  author  = {Alegr{\'i}a, Alfredo and Emery, Xavier and Filosi, Tobia and Porcu, Emilio},
  title   = {Computationally efficient algorithms for simulating isotropic {G}aussian random fields on graphs with {E}uclidean edges},
  journal = {J. Comp. Graph. Statist.},
  year    = {2026},
  doi     = {10.1080/10618600.2025.2574535},
  note    = {Advance online publication},
}

@article{BSW2022,
  author  = {Bolin, David and Simas, Alexandre B. and Wallin, Jonas},
  title   = {Gaussian {W}hittle--{M}at\'ern fields on metric graphs},
  journal = {Bernoulli},
  volume  = {30},
  number  = {2},
  pages   = {1611--1639},
  year    = {2024},
  doi     = {10.3150/23-BEJ1647}
}

@article{BSW2023_AOS,
  author  = {Bolin, David and Simas, Alexandre B. and Wallin, Jonas},
  title   = {Statistical inference for Gaussian {W}hittle--{M}at\'ern fields on metric graphs},
  journal = {J. R. Statist. Soc. B},
  year    = {2026},
  doi     = {10.1093/jrsssb/qkag074},
  note    = {Advance article}
}

@article{BSW_Markov,
  author  = {Bolin, David and Simas, Alexandre B. and Wallin, Jonas},
  title   = {Markov properties of Gaussian random fields on compact metric graphs},
  journal = {Bernoulli},
  volume  = {32},
  number  = {1},
  pages   = {153--178},
  year    = {2026},
  doi     = {10.3150/25-BEJ1853}
}

@article{anderes2020isotropic,
  author  = {Anderes, Ethan and M{\o}ller, Jesper and Rasmussen, Jakob G.},
  title   = {Isotropic covariance functions on graphs and their edges},
  journal = {Ann. Statist.},
  volume  = {48},
  number  = {4},
  pages   = {2478--2503},
  year    = {2020},
  doi     = {10.1214/19-AOS1896}
}

@article{Hoef2006,
  author  = {Ver Hoef, Jay M. and Peterson, Erin E. and Theobald, David},
  title   = {Spatial statistical models that use flow and stream distance},
  journal = {Environ. Ecol. Statist.},
  volume  = {13},
  number  = {4},
  pages   = {449--464},
  year    = {2006},
  doi     = {10.1007/s10651-006-0022-8}
}

@article{moller2022lgcp,
  author  = {M{\o}ller, Jesper and Rasmussen, Jakob G.},
  title   = {Cox processes driven by transformed Gaussian processes on linear networks---A review and new contributions},
  journal = {Scand. J. Statist.},
  volume  = {51},
  number  = {3},
  pages   = {1288--1322},
  year    = {2024},
  doi     = {10.1111/sjos.12720}
}

@article{AliniaeifardWangVanWilligenburg2024,
  author  = {Farid Aliniaeifard and Victor Wang and Stephanie van Willigenburg},
  title   = {Deletion-contraction for a unified Laplacian and applications},
  journal = {Linear Algebra Appl.},
  volume  = {691},
  pages   = {50--95},
  year    = {2024},
  doi     = {10.1016/j.laa.2024.03.013}
}

@inproceedings{borovitskiy2021,
  author    = {Borovitskiy, Viacheslav and Azangulov, Iskander and Terenin, Alexander and Mostowsky, Peter and Deisenroth, Marc and Durrande, Nicolas},
  title     = {Mat{\'e}rn {G}aussian processes on graphs},
  booktitle = {Proc. 24th Int. Conf. Artif. Intel. Statist.},
  series    = {Proc. Mach. Learn. Res.},
  volume    = {130},
  pages     = {2593--2601},
  publisher = {PMLR},
  year      = {2021},
}

@incollection{kostenko_nicolussi2023,
  author    = {Kostenko, Aleksey and Nicolussi, Noema},
  title     = {Laplacians on infinite graphs: discrete vs. continuous},
  booktitle = {European Congress of Mathematics},
  publisher = {EMS Press},
  pages     = {295--323},
  year      = {2023},
  doi       = {10.4171/8ECM/02}
}

@article{gneiting_raftery2007,
  author = {Gneiting, Tilmann and Raftery, Adrian E.},
  title = {Strictly proper scoring rules, prediction, and estimation},
  journal = {J. Am. Statist. Assoc.},
  volume = {102},
  number = {477},
  pages = {359--378},
  year = {2007}
}

@misc{paris_comptages,
  author = {{Ville de Paris}},
  title = {Comptage routier -- {D}onn\'ees trafic issues des capteurs permanents},
  howpublished = {\url{https://opendata.paris.fr/explore/dataset/comptages-routiers-permanents}},
  note = {Accessed 23 October 2026},
  year = {2026}
}

@misc{webtris,
  author = {{National Highways}},
  title = {{WebTRIS}: traffic flow and speed data for the strategic road network},
  howpublished = {\url{https://webtris.nationalhighways.co.uk}},
  note = {Accessed 22 October 2026},
  year = {2026}
}

@misc{openstreetmap,
  author = {{OpenStreetMap contributors}},
  title = {{OpenStreetMap}},
  howpublished = {\url{https://www.openstreetmap.org}},
  note = {Data retrieved through the Overpass API on 22 October 2026},
  year = {2026}
}

@article{filosi2026polyharmonic,
  title = {Smooth Isotropic Covariance Functions on Metric Graphs via Polyharmonic Resistance Distances},
  author = {Filosi, Tobia and Porcu, Emilio and Agostinelli, Claudio},
  year = {2026},
  journal = {arXiv preprint arXiv:2609.33519},
  doi = {10.48550/arXiv.2609.33519}
}

\end{document}